\documentclass[12pt,a4paper]{article}

\usepackage{amsmath,amssymb}
\usepackage{xcolor}

\def\cD{{\cal D}}

\def\cA{{\cal A}}

\def\CP{{{\mathbb C}{\rm P}}}
\def\RP{{{\mathbb R}{\rm P}}}
\def\Aut{{\rm Aut}}
\def\bdeg{{\rm bdeg}}

\def\R{{\mathbb R}}
\def\C{{\mathbb C}}
\def\Q{{\mathbb Q}}

\def\p{{\partial}}

\usepackage{theorem}

\newtheorem{theorem}{Theorem}[section]
\newtheorem{proposition}{Proposition}[section]
\newtheorem{corollary}[proposition]{Corollary}

{\theorembodyfont{\rmfamily}
\newtheorem{definition}[proposition]{Definition}

\newtheorem{remark}[proposition]{Remark}

}

\usepackage{graphicx}

\title{On framed simple purely real Hurwitz numbers}

\author{
Maxim Kazarian\thanks{Steklov Mathematical Institute RAS,
National Research University Higher School of Economics,
Skolkovo Institute of Science and Technology, kazarian@mccme.ru},
Sergey Lando\thanks{National Research University Higher School of Economics,
Skolkovo Institute of Science and Technology, lando@hse.ru},
Sergey Natanzon\thanks{National Research University Higher School of Economics
and Institute of Theoretical and Experimental Physics, natanzons@mail.ru}}

\date{}

\begin{document}
\maketitle

\begin{abstract}
We present a study of real Hurwitz numbers enumerating a special kind
of real meromorphic functions, which we call simple framed purely real functions.
We deduce partial differential equations of
cut-and-join type for generating functions for these numbers.
We also construct a topological field theory for them.
\end{abstract}

A real algebraic curve is a complex curve together with an antiholomorphic involution.
Points invariant with respect to the involution are the real points of the curve.
A meromorphic function $f:C\to\CP^1$ on a real algebraic curve~$C$ is real provided it is equivariant
with respect to the antiholomorphic involution on~$C$ and complex conjugation on the target.
We say that a real function~$f$ is simple if all its finite critical values are
simple.

In the complex case, Hurwitz numbers enumerate meromorphic functions with a given set of critical values,
ramification over each being a prescribed partition of the degree of the function.
Hurwitz numbers do not depend on the specific positions of the critical values.
Simple Hurwitz numbers enumerate meromorphic functions with a given set of critical values,
ramification over one of which is a prescribed partition of the degree of the function, while
all the other critical values are simple. Simple Hurwitz numbers play a crucial role in the
study of intersection theory on moduli spaces of algebraic curves. We usually refer to the
point in~$\CP^1$ over which the ramification data is prescribed as to infinity, and its preimages are
poles. The parts of the partition are, therefore, the orders of the poles.
All the other critical values are finite.

In contrast to the complex case, the numbers of real meromorphic functions
with prescribed partitions over critical values
depend essentially on the location of the critical values.
This is not true, however, for real simple Hurwitz numbers, which
enumerate real meromorphic functions
with prescribed orders of the poles and a given set of real finite critical values,
all being simple.
A classical example
of such numbers are Bernoulli and Euler numbers (that is, numbers of up-down permutations),
which enumerate simple real polynomials of given degree.
The goal of this paper is to extend our understanding of up-down permutations
to real functions on separating real curves of arbitrary genus.
We hope that our results could be helpful in understanding the geometry of moduli space of real curves
with marked points.

When working on this paper, S.N. was partly supported by the RFBR Grant 16-01-00409
and enjoyed the hospitality of Max-Planck Institut f\"ur Mathematik (Bonn) and Institut
des Hautes \'Etude Scientifique (Paris).

\section{Definitions and statement of the main result}\label{sDS}
\subsection{Simple separating real meromorphic functions}

Let~$C$ be a compact nonsingular complex algebraic curve. For an anti-holomorphic
involution $\tau:C\to C$, the pair $(C,\tau)$ is called a (compact nonsingular)
{\it real algebraic curve}. Of course, an anti-holomorphic involution exists
not on every complex curve; on the other hand, some complex curves admit several
nonequivalent anti-holomorphic involutions, whence several real curve structures.
A detailed study of multiple real structures on complex curves can be found in~\cite{N88,N04}.

The fixed points of the involution~$\tau$
form the {\it real part} $C^\tau$ of the curve~$(C,\tau)$, $C^\tau\subset C$.
The real part $C^\tau$ is a $1$-dimensional real manifold, whence
a disjoint union of circles. A connected real curve $(C,\tau)$ is said to be {\it separating\/} provided
$C\setminus C^\tau$ is disconnected, and {\it nonseparating\/} otherwise. For a separating real curve,
 the complement $C\setminus C^\tau$ consists of two connected components interchanged by the involution~$\tau$. By a {\it framing\/}
 of a separating real curve we mean a choice of one of the two components of $C\setminus C^\tau$.
 A framing chosen, the real part $C^\tau$
 of a separating real curve, being a $1$-dimensional real manifold, acquires a natural orientation
 as the boundary of the chosen component of its complement equipped with the complex orientation. A change of the framing of a separating real curve
 leads to simultaneous change of orientations of all connected components of its real part.

A {\it real holomorphic mapping\/} of a real curve $(C_1,\tau_1)$ to a real curve $(C_2,\tau_2)$
is a holomorphic mapping $f:C_1\to C_2$ equivariant with respect to the pair of involutions $(\tau_1,\tau_2)$,
that is, such that $f\circ\tau_1=\tau_2\circ f$. In particular, a {\it real meromorphic function\/} on a
real curve $(C,\tau)$ is a real holomorphic mapping from $(C,\tau)$ to $(\CP^1,\sigma)$,
where $\sigma:\CP^1\to\CP^1$ is the standard complex conjugation.

A real meromorphic function $f:C\to\CP^1$ on a real curve $(C,\tau)$
is said to be {\it simple\/} if all its finite critical values
are simple.
A real meromorphic function $f:C\to\CP^1$ on a real curve $(C,\tau)$
is said to be {\it purely real\/} if all its finite critical values
are real.
The main object of our study in this paper are {\it simple purely real meromorphic
functions on framed separating real curves}.

%\medskip
%{\it Ramification type of a critical value. To each critical value we associate its ramification type~$\mu$ which can be encoded in many different ways: collection of partitions, partitions with additional markings etc. It should be introduced in such a way that the notation could be easily adjusted when we pass from oriented to usual functions. A possible way to do it is to associate a monomial $p_\mu$ in the variables $p_1^\pm$, $p_2^\pm$, \dots, $q_1$, $q_2$, \dots. The degree of the monomial is equal to the degree of the function where $\deg p_k^+=\deg p_k^-=k$, $\deg q_k=2k$.}

%\emph{To each ramification type $\mu$ of a pure real function we associate a monomial $p_\mu$. For the case of a general real surface $p_\mu$ is a monomial in the variables $p_{2k}^\pm$, $p_{2k+1}$, $q_k$. For the case of a function on an oriented real surface it is a monomial in the variables $p_k^\pm$ and $q_k$. In both cases the degree of the monomial is equal to $|\mu|$ with respect to the grading with $\deg p_k^\pm=1$, $\deg q_k=2$.}

%\emph{Ramification type of a trivial ramification and an associated monomial}

%\emph{Ramification type of a simple ramification and an associated monomial}

\subsection{Framed purely real simple Hurwitz numbers}\label{sFSHN}

By a {\it framed real meromorphic function\/} we understand a real meromorphic
function $f:(C,\tau)\to(\CP^1,\sigma)$ defined on a framed separating real curve~$(C,\tau)$.
We denote the connected component of the complement $C\setminus C^\tau$ chosen by the
framing by~$C^f$.
Two framed real meromorphic functions $f_1:C_1\to\CP^1$, $f_2:C_2\to\CP^1$ on framed separating real holomorphic curves
$(C_1,\tau_1)$, $(C_2,\tau_2)$ are said to be {\it equivalent\/} if there is an invertable real holomorphic
mapping $\varphi:C_1\to C_2$ such that $f_1=f_2\circ \varphi$ and $\varphi(C_1^{f_1})=C_2^{f_2}$.
In particular, the domains of two equivalent meromorphic functions have the same genus.
For a given degree of a function, a given genus of its domain, and a given set
of critical values, the set of equivalence classes of functions possessing this data is finite.

Let~$f:(C,\tau)\to(\CP^1,\sigma)$ be a framed real meromorphic function;
we are going to define its ramification type over a point in~$\RP^1$.
Without loss of generality we may suppose that this point is~$\infty\in\RP^1$,
so that its preimages are the poles of~$f$.

The poles of~$f$ are split into real ones
and pairs of $\tau$-conjugate nonreal poles. In each pair of $\tau$-conjugate
poles, exactly one of them belongs to the domain~$C^f$ of the function~$f$.
The orders of the $\tau$-conjugate poles in the domain~$C^f$ form a
partition~$\lambda=(\ell_1,\ell_2,\dots)$.
Real poles of~$f$ are split into two types, which we call positive and negative.
A real pole is said to be {\em positive\/} (respectively, {\em negative})
if the function~$f$ increases (respectively, decreases) to the left of the pole
(with respect to the orientation of~$C^\tau$). Note that knowing the type of a pole
and the parity of its order one can reconstruct the behavior of~$f$ to the right of the pole
(whether~$f$ is increasing or decreasing).
The orders of positive (respectively, negative) real poles of~$f$ form a partition
$\kappa^+=(k^+_1,k^+_2,\dots)$ (respectively, $\kappa^-=(k^-_1,k^-_2,\dots)$).
Together, the partitions $\kappa^+,\kappa^-$, and $2\lambda=(2\ell_1,2\ell_2,\dots)$ form a partition
of the degree of~$f$. Below, we use the multiplicative form of writing out partitions,
so that $1^{a_1}2^{a_2}\dots$ denotes a partition with~$a_1$ parts equal to~$1$,
$a_2$ parts equal to~$2$, and so on.

The triple of partitions $\mu=(\kappa^+,\kappa^-,\lambda)$
is called the {\em ramification type\/} of~$f$ over infinity.
If all the poles are of order~$1$,
then each of the three partitions $\kappa^+,\kappa^-,\lambda$ consists
of parts~$1$.

We are interested in the case when all finite critical values are real and simple.
Denote by $h^\R_{m;\mu}$ the {\em framed simple purely real connected Hurwitz number\/}
enumerating purely real meromorphic functions with connected domain having
ramification type $\mu=(\kappa^+,\kappa^-,\lambda)$ over infinity and~$m$ given nondegenerate real critical values.
Formally,
$$
h^\R_{m;\mu}=\sum_{[f]}\frac1{\#\Aut[f]},
$$
where the summation is carried over all equivalence classes $[f]$ of framed simple purely real meromorphic functions~$f$
defined on connected curves and having a prescribed set of~$m$ finite nondegenerate
real critical values and ramification of type~$\mu$ over infinity;
here $\#\Aut[f]$ denotes the order of the automorphism group of the equivalence class.
We denote by $h^{\R\circ}_{m;\mu}$
the number of all simple framed purely real meromorphic functions, that is, including those whose domain is not necessarily connected.

Note, however, that the automorphism group of a framed purely real function with a {\em connected\/} domain
is trivial, so that the order of its automorphism group is~$1$.
For example, if $C$ is a genus~$g\ge2$ hyperelliptic curve,
which is the compactification of the curve
$$
y^2=P_{2g+1}(x),
$$
in~$\C^2$, where~$P_{2g+1}(x)$ is a real polynomial of degree~$2g+1$ with $2g+1$ pairwise distinct real roots,
then the function on this curve taking a point to its~$x$ coordinate admits a nontrivial automorphism,
which is the hyperelliptic involution $(x,y)\mapsto(x,-y)$.
This automorphism, however, exchanges the two framings of the function, so that
$$
h^\R_{2g-1;2^1,\emptyset,\emptyset}=h^\R_{2g-1;\emptyset,2^1,\emptyset}=1.
$$
A framed purely real function~$f$ with a disconnected domain can have nontrivial automorphisms only if
the domain~$C^f$ contains several connected components of genus~$0$ without boundary, the restriction of~$f$ to each
of these components being of degree~$1$. In this case, the automorphism group of~$f$
consists of permutations of such connected components with a positive pole, and, separately, with
a negative pole.

It is convenient to make an additional convention by setting
$$h^{\R}_{0;\emptyset,\emptyset,1}=1.$$
This Hurwitz number corresponds to the degree two covering of $\C P^1$ such that the covering surface consists of two disjoint components both of which are mapped to the target Riemann sphere isomorphically, and the involution~$\tau$ interchanges these two components.
This convention implies, in particular,
 that this covering should be considered as a real function with a \emph{connected} domain of Euler characteristic~$4$ with no real points.
 As a corollary, when considering functions on disconnected curves we allow them
 to have several copies of the above mapping.
 The reason for this convention is in a simpler form of the differential equation for the generating functions. A
  deeper reason to follow it will be explained in Sec.~\ref{sTFT}.

Let us associate with the ramification type $\mu=(\kappa^+,\kappa^-,\lambda)=((k_1^+,k_2^+,\dots),(k_1^-,k_2^-,\dots),(\ell_1,\ell_2,\dots))$ the monomial
  $$ p_\mu=p_{k_1^+}^+p_{k_2^+}^+\dots p_{k_1^-}^-p_{k_2^-}^-\dots q_{\ell_1}q_{\ell_2}\dots$$
in the independent commuting variables $p_k^+,p_k^-,q_k$, $k=1,2,\dots$. Introduce the generating functions
\begin{eqnarray*}
H^{\R}(u,\hbar;p^+_1,\dots,p^-_1,\dots,q_1,\dots)&=&\sum\hbar^{\chi_{m;\mu}} h^\R_{m;\mu}p_\mu\frac{u^m}{m!},\\
H^{\R\circ}(u,\hbar;p^+_1,\dots,p^-_1,\dots,q_1,\dots)&=&\sum\hbar^{\chi_{m;\mu}} h^{\R\circ}_{m;\mu}p_\mu\frac{u^m}{m!},
\end{eqnarray*}
where the summation on the right runs over all triples of partitions $\mu=(\kappa^+,\kappa^-,\lambda)$
and all nonnegative values of~$m$. Here $\chi_{{m;\mu}}$ denotes the Euler characteristics of the
source curve~$C$ of the meromorphic function~$f$.

\begin{remark}\label{remh}
Introduction of the explicit parameter~$\hbar$ whose power indicates the Euler characteristic~$\chi$
of the source of the function is not the only way to determine this Euler characteristic from
the Hurwitz numbers. Indeed, the Riemann--Hurwitz formula states that for a function with the
pole partition~$\mu=(\kappa^+,\kappa^-,\lambda)$ and $m$ simple finite critical points
$$
\chi_{{m;\mu}}= |\kappa^++\kappa^-+2\lambda|+\ell(\kappa^+)+\ell(\kappa^-)+2\ell(\lambda)-m,
$$
where $|\kappa^++\kappa^-+2\lambda|$ is the degree of~$f$, and $\ell(\cdot)$ denotes the length of a partition, respectively. In other words, the functions $H^{\R}$, $H^{\R\circ}$ can be recovered from their specializations at $\hbar=1$ by the rescaling of the variables $p_k^\pm\mapsto\hbar^{k+1}p_k^\pm$, $q_k\mapsto\hbar^{2k+2}q_k$, $u\mapsto \hbar^{-1}u$. However, it is more convenient to follow an explicit indication of the Euler characteristic
in the generating function.
\end{remark}

As usual, we have the following relationship between the generating functions enumerating connected
and disconnected meromorphic functions:
$$
H^{\R\circ}=\exp\left( H^\R\right).
$$

\begin{theorem}\label{tcjdpm}
The generating function $H^{\R\circ}$ satisfies the differential equation
\begin{equation}\label{ecjdpm}
\frac{\p H^{\R\circ}}{\partial u}=W^+(H^{\R\circ});
\end{equation}
here the linear partial differential operator~$W^+$ is defined as follows:
$$
W^+=\sum_{i,j=1}^\infty\left(p_i^{{\bar i}}p^+_j\frac\p{\p p^{{\bar i}}_{i+j}}
+\hbar^{-2} p^{{\bar i}}_{i+j}\frac{\p^2}{\p p_i^{{\bar i}}\p p^+_j}\right)
+\sum_{i=1}^\infty\left(i\hbar^{-2} p^+_{2i}\frac{\p}{\p q_i}+q_i\frac{\p}{\p p^+_{2i}}\right),
$$
where, for a positive integer~$i$, notation ${\bar i}$ stands for the sign~$+$ provided~$i$ is even, and
for the sign~$-$ otherwise.
\end{theorem}

If we consider the generating function~$H^{\R\circ}$ as a power
series in~$u$, then Eq.~(\ref{ecjdpm}) together with the initial conditions at $u=0$, which are
$$
H^{\R\circ}(0,\hbar;p_1^\pm,p_2^\pm,\dots)=e^{\hbar^{-2}(p^+_1+p^-_1)+\hbar^{-4}q_1},
$$
allow one to compute as many terms of the power
series expansions as we like. Indeed, we have
\begin{eqnarray*}
W^+\left(e^{\hbar^{-2}(p^+_1+p^-_1)+\hbar^{-4}q_1}\right)&=&e^{\hbar^{2}(p^+_1+p^-_1)+\hbar^{-4}q_1}\hbar^{-2}(p^+_2+p^-_2)\\
W^+\left(e^{\hbar^{-2}(p^+_1+p^-_1)+\hbar^{-4}q_1}\hbar^{-2}(p^+_2+p^-_2)\right)&=&e^{\hbar^{-2}(p^+_1+p^-_1)+\hbar^{-4}q_1}
\hbar^{-4}
\left((p^+_2+p^-_2)^2\right.\\
&&\left.+\hbar^{2}(p^+_3+p^-_3+p^+_1p^-_1+q_1)\right),\\
\end{eqnarray*}
and so on, so that
\begin{eqnarray*}
H^{\R\circ}(u,\hbar;p_1^\pm,p_2^\pm,\dots)&=&e^{\hbar^{-2}(p^+_1+p^-_1)+\hbar^{-4}q_1}\left(1+\hbar^{-2}(p^+_2+p^-_2)\frac{u}{1!}\right.\\
&&\left.+\hbar^{-4}\left((p^+_2+p^-_2)^2+\hbar^{2}(p^+_3+p^-_3+p^+_1p^-_1+q_1)\right)\frac{u^2}{2!}+\dots\right).
\end{eqnarray*}

This recursive procedure can also be written in the closed form
$$
H^{\R\circ}(u,\hbar;p_1^\pm,p_2^\pm,\dots)=e^{u\,W^+}e^{\hbar^{-2}(p^+_1+p^-_1)+\hbar^{-4}q_1}.
$$

The leading terms of the generating function~$H^{\R\circ}$ immediately give, by taking the
logarithm, the leading terms of the generating function~$H^\R$:
\begin{eqnarray*}
H^{\R}(u,\hbar;p_1^\pm,p_2^\pm,\dots)&=&\hbar^{-4}q_1+\hbar^{-2}(p^+_1+p^-_1)+\hbar^{-2}(p^+_2+p^-_2)\frac{u}{1!}\\
&&+\hbar^{-2}\left(p^+_3+p^-_3+p^+_1p^-_1+q_1\right)\frac{u^2}{2!}\\
&&+\hbar^{-2}\left((p^+_2+p^-_2)(p^+_1+p^-_1)+2(p^+_4+p^-_4)+\hbar^{-2}(p^+_2+p^-_2)\right)\frac{u^3}{3!}+\dots.
\end{eqnarray*}

The logarithm $H^\R$ of the generating function $H^{\R\circ}$ satisfies a partial differential
equation, which can be deduced by substituting $H^{\R\circ}=e^{H^\R}$
into the equation in the theorem. The partial differential equation for~$H^\R$
is linear no longer, but we can use it to deduce the power series expansion for~$H^\R$
as well.

\begin{corollary}
The generating function $H^{\R}$ satisfies the differential equation
\begin{eqnarray*}\label{ecjdpmc}
\frac{\partial H^{\R}}{\partial u}&=&\sum_{i,j=1}^\infty\left(p_i^{{\bar i}}p^+_j\frac{\p H^\R}{\p p^{{\bar i}}_{i+j}}
+p_{i+j}^{{\bar i}}\left(\hbar^{-2} \frac{\p^2 H^\R}{\p p_i^{{\bar i}}\p p^+_j}+\frac{\p H^\R}{\p p_i^{{\bar i}}}\frac{\p H^\R}{\p p_j^+}\right)\right)\\
&&+\sum_{i=1}^\infty\left(i\hbar^{-2} p^+_{2i}\frac{\p H^\R}{\p q_i}+q_i\frac{\p H^\R}{\p p^+_{2i}}\right),
\end{eqnarray*}
where, for a positive integer~$i$, notation ${\bar i}$ stands for~$+$ provided~$i$ is even, and for~$-$
otherwise.
\end{corollary}

Note that in spite of the fact that the operator~$W^+$ is not symmetric with respect to the
exchange of the variables~$p^+_k\leftrightarrow p^-_k$ for $k=1,2,\dots$, the generating function~$H^{\R\circ}$
(and hence the generating function~$H^\R$) possesses this property.

\begin{corollary}\label{cconn}
The generating function $H^{\R\circ}$ satisfies the differential equations
$$
\frac{\p H^{\R\circ}}{\partial u}=W^-(H^{\R\circ}),\qquad
\frac{\p H^{\R\circ}}{\partial u}=W(H^{\R\circ}),
$$
where the partial differential operator $W^-$ is obtained from $W^+$ by replacing
each variable~$p^+_k$ by~$p^-_k$ and vice versa, for~$k=1,2,\dots$, and where
$$
W=\frac12(W^++W^-).
$$
\end{corollary}

The operator~$W$ {\em is\/} symmetric with respect to the exchange of the variables $p_k^+\leftrightarrow p_k^-$
but it contains more terms than~$W^+$ and is therefore less efficient from the practical point of view.
It is easy to verify that the operators $W^+$ and~$W^-$ commute with one another, whence
both of them commute with~$W$.

Equation~(\ref{ecjdpm}) is a cut-and-join type equation similar to the one in~\cite{GJ97} for
the generating function for ordinary, that is complex, simple Hurwitz numbers. The equation
in~\cite{GJ97} is simpler due to the fact that there is no difference between the
real and the complex poles, as well as between positive and negative poles,
which allows one to have in the complex case a single infinite sequence of
variables instead of three sequences in the real case.

Similarly to the complex cut-and-join equation, the differential operator~$W^+$ on the right (as well as $W^-$ and~$W$)
is homogenous, meaning that it preserves the subspaces of polynomials of given homogeneous degree.
In our real case the degree of a monomial is defined as the sum of the degrees of the variables
it contains, which are $\deg p^\pm_{k}=k$, $\deg q_k=2k$, $k=1,2,\dots$. The degree of a monomial coincides
with the degree of a meromorphic function contributing to its coefficient.
%As in the complex case,
%the operators~$W^-,W^+,W$ are self-adjoint with respect to a nondegenerate scalar product
%on the space of polynomials of given degree, see Sec.~\ref{sTFT}.

In the real case, however, each of the operators~$W^-,W^+,W$ preserves additionally a finer
splitting of the space of polynomials in the variables $p^\pm_k,q_k$ given by the following bigrading: set
\begin{equation}
\bdeg~ p_{2k}^\pm=\bdeg~ q_k=(k,k),\quad\bdeg~ p_{2k+1}^+=(k+1,k),\quad\bdeg~ p_{2k+1}^-=(k,k+1).
\end{equation}\label{ebd}
This statement can be easily verified since each summand in the definition of the operator~$W^+$
(whence of both other operators) preserves the indicated bidegree.
The bidegree matches the natural bidegree of a framed separating real function: for a given such function
$f:(C,\tau)\to(\CP^1,\sigma)$, its bidegree is the pair (the number of preimages of the upper hemisphere
in~$C^f$, the number of preimages of the upper hemisphere
in the complement~$C\setminus C^f$).
Since the operator $W^+$ of Theorem~\ref{tcjdpm}
 preserves the bigrading, it acts  in the space of polynomials of fixed bidegree.
  In other words, the differential equation~\eqref{ecjdpm} splits as a direct sum of linear ordinary differential equations with constant coefficients in the finite dimensional spaces of polynomials of fixed bidegree, and the exponential of the operator $W^+$ can be computed as the direct sum of exponentials of its restriction to the subspaces of fixed bidegree.
  In Sec.~\ref{sTFT} we introduce a nondegenerate scalar product on each subspace of polynomials of fixed bidegree,
  and show that the restrictions of the operators~$W^+,W^-,W$ to these subspaces are self-adjoint,
  whence diagonalizable, operators.

\begin{remark}
Up to now, our attempts to deduce simpler evolution equations for functions obtained
from~$H^{\R\circ}$ by getting rid of complex poles or of the distinction between the two
types of real poles failed. In the latter case, we obtain a simpler equation only for genus~$0$,
see Sec.~\ref{sg0}.
\end{remark}

\begin{remark}
A close problem of computing disc simple Hurwitz numbers was treated in~\cite{N10} by similar methods.
Moreover, in addition to real framed meromorphic functions,~\cite{N10} considers more general Dold--Smith coverings corresponding to real meromorphic functions on not necessary separating real curves, see Appendix~A.
In this case both the differential equations and the initial conditions become more complicated, which makes
presumable computations less efficient.
\end{remark}

\section{Genus~$0$ case}\label{sg0}

In this section we analyze specification of the generating function of framed simple real Hurwitz numbers
to the case of rational functions $(\CP^1,\sigma)\to(\CP^1,\sigma)$.

\subsection{Cut-and-join for rational functions}
The generating function for framed connected simple Hurwitz numbers with the genus~$0$ domain,
that is, the case of rational functions, is the coefficient of~$\hbar^2$ in~$H^\R$.
Denote by~$H_0^\R$ the result of the substitution $p^+_i=p_i$, $p^-_i=p_i$, for $i=1,2,\dots$
into one half of this coefficient, so that
$$
H_0^\R(u;p_1,p_2,\dots,q_1,q_2,\dots)=\frac12[\hbar^2]H^\R(u,\hbar;p_1,p_2,\dots,p_1,p_2,\dots,q_1,q_2,\dots).
$$
In the generating function~$H_0^\R$, we make no difference between the positive and the negative
real poles.

In the genus~$0$ case the cut-and-join equations simplify a lot since they do not involve second order partial derivatives.
Corollary~\ref{cconn} immediately implies
\begin{theorem}
The function $H_0^\R$ satisfies the partial differential equation
$$
\frac{\p H_0^{\R}}{\p u}=
\frac12\sum_{i,j=1}^\infty\left( p_ip_j\frac{\p H_0^{\R}}{\p p_{i+j}}
+p_{i+j}\frac{\p H_0^{\R}}{\p p_{i}}\frac{\p H_0^{\R}}{\p p_j}\right)
+\frac12\sum_{i=1}^\infty q_i\frac{\p H_0^\R}{\p p_{2i}}+\frac12p_2.
$$
\end{theorem}

\subsection{Comparison with known enumerative results}

There are very few enumerative results concerning real Hurwitz numbers,
all of them related to the genus~$0$ case.

Simple real polynomials of degree~$n$ are enumerated by Bernoulli (for~$n$ odd) and Euler (for~$n$ even) numbers,
so that the exponential generating function has the form
$$
\frac1{\cos x}+\tan x=1+\frac{x}{1!}+\frac{x^2}{2!}+2\frac{x^3}{3!}+5\frac{x^4}{4!}+16\frac{x^5}{5!}+61\frac{x^6}{6!}+272\frac{x^7}{7!}+\dots.
$$
These numbers $1,1,1,2,5,16,61,\dots$ appear as coefficients of the monomials $p^\pm_n\frac{u^{n-1}}{(n-1)!}$ both in the generating functions~$H^\R$
and $H^{\R\circ}$, and as coefficients of the monomials $p_n\frac{u^{n-1}}{(n-1)!}$ in~$H_0^\R$.

Real generic rational functions of degree~$n$ were counted in~\cite{SV03}.
A meromorphic function is said to be {\em generic\/} if all its poles have order~$1$.
The sequence enumerating generic real functions is computed in~\cite{SV03}
and, starting with $n=3$, begins with the
numbers
$$
2,20,406,\dots.
$$
These numbers are obtained by summing up the coefficients
of the monomials $p_1^{n-2k}q_1^k \frac{u^{2n-2}}{(2n-2)!}$, $k=0,1,\dots,[n/2]$
in the generating function~$H_0^\R$.
The corresponding generating function has the form
$$
H_0^\R(u;1,0,0,\dots,1,0,0,\dots)=1+\frac12\frac{u^2}{2!}+2\frac{u^4}{4!}+20\frac{u^6}{6!}+
406\frac{u^8}{8!}+14652\frac{u^{10}}{10!}+\dots.
$$

\section{Diagrammatic description of real meromorphic functions}\label{sDD}

Our proof of the main Theorem~\ref{tcjdpm} is based on the diagrammatic
description of real meromorphic functions
originating in~\cite{NSV02}. This approach was developed in~\cite{CNS16}.
Note that a similar diagrammatic approach was applied by S.~Barannikov in~\cite{B92}
in the specific situation of polynomials.
Barannikov established a topological classification of real polynomials all whose
critical values (both real and nonreal ones) are simple.

Since the terminology is not yet stable (the authors use such notions as
``chord diagrams'', ``gardens'', ``parks'', and so on) we take liberty
to suggest just the term ``diagram of a framed real meromorphic function''
for our current purposes. (In particular, the first two authors of the present paper are used to
apply the term ``chord diagram'' in a completely different environment, namely,
in Vassiliev's theory of finite order knot invariants).

\subsection{Abstract diagrams}

Abstract diagrams
are aimed at being in one-to-one correspondence with topological types of framed real
meromorphic functions. We start with introducing the notion of an abstract diagram.

\begin{definition}\label{defdiag}
An {\em abstract diagram\/} is a compact oriented two-dimensional surface $C_+$ with boundary
such that all its connected components have nonempty boundary, together with an oriented graph drawn on it
possessing the following properties:
\begin{itemize}
\item the vertices of the graph belonging to the boundary are split into two different classes,
called \emph{critical points} and \emph{poles}, respectively; all the other vertices are called {\em internal poles};
\item each interval between two
consecutive vertices on the boundary is an edge of the graph (such edges are said to be {\em boundary};
all the other edges are said to be {\em internal});
\item there is exactly one internal edge incident to each critical point;
\item the complement to the graph on the surface is a disjoint union of open disks; the boundary of
each disk passes through a pole (either boundary or internal one) exactly once;
below, we call these discs the {\em faces\/} of the diagram;
\item the graph is oriented in such a way that the orientation of the edges along
the boundary of each face is consistent;
\item let~$m$ denote the number of critical points; the critical points are numbered from~$1$ to~$m$
in such a way that the numbering increases along each oriented edge connecting two critical points.
\end{itemize}
\end{definition}

Note that orientation of the edges of a diagram has nothing to do with the canonical
orientation of the boundary, and normally each connected component of the boundary
of a diagram is split into several boundary edges oriented in an alternating order.
The faces of any diagram admit a chessboard coloring: we may color white those
faces the orientation of whose boundaries coincides with the one induced by the orientations
of the faces, and color black those faces the orientation of whose boundaries is opposite
to the one induced by the orientations of the faces.

By the {\em genus\/} of a connected abstract diagram we mean the genus of the
compact oriented surface without boundary obtained by taking the double
of the underlying surface of the diagram.

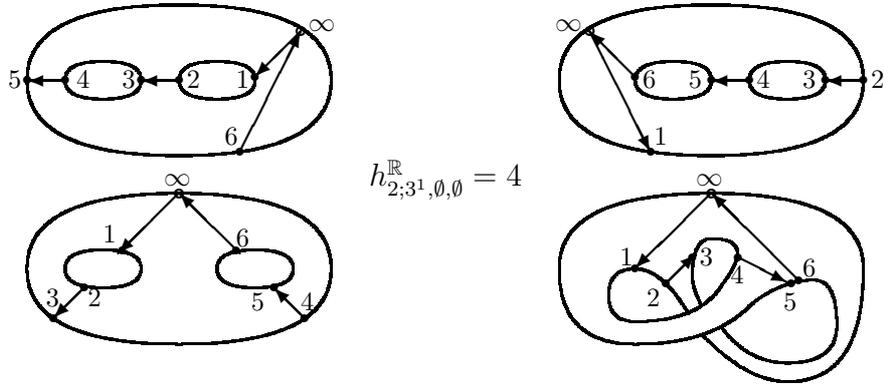
\begin{figure}
\setlength{\unitlength}{0.5mm}
\thicklines
\begin{picture}(400,100)(-25,0)
\qbezier(0,80)(0,100)(40,100)
\qbezier(80,80)(80,100)(40,100)
\qbezier(0,80)(0,60)(40,60)
\qbezier(80,80)(80,60)(40,60)
\qbezier(10,80)(10,85)(20,85)
\qbezier(30,80)(30,85)(20,85)
\qbezier(10,80)(10,75)(20,75)
\qbezier(30,80)(30,75)(20,75)
\qbezier(40,80)(40,85)(50,85)
\qbezier(60,80)(60,85)(50,85)
\qbezier(40,80)(40,75)(50,75)
\qbezier(60,80)(60,75)(50,75)
\put(10,80){\vector(-1,0){10}}
\put(40,80){\vector(-1,0){10}}
\put(72,93){\vector(-1,-1){12}}
\put(56,61){\vector(1,2){16}}
\put(72,93){\circle{2}}
\put(56,61){\circle*{2}}
\put(10,80){\circle*{2}}
\put(0,80){\circle*{2}}
\put(30,80){\circle*{2}}
\put(40,80){\circle*{2}}
\put(60,81){\circle*{2}}
\put(74,93){\mbox{\footnotesize $\infty$}}
\put(52,63){\mbox{\footnotesize $6$}}
\put(13,78){\mbox{\footnotesize $4$}}
\put(-5,78){\mbox{\footnotesize $5$}}
\put(25,78){\mbox{\footnotesize $3$}}
\put(42,78){\mbox{\footnotesize $2$}}
\put(55,78){\mbox{\footnotesize $1$}}
\qbezier(140,80)(140,100)(180,100)
\qbezier(220,80)(220,100)(180,100)
\qbezier(140,80)(140,60)(180,60)
\qbezier(220,80)(220,60)(180,60)
\qbezier(210,80)(210,85)(200,85)
\qbezier(200,85)(190,85)(190,80)
\qbezier(210,80)(210,75)(200,75)
\qbezier(200,75)(190,75)(190,80)
\qbezier(180,80)(180,85)(170,85)
\qbezier(170,85)(160,85)(160,80)
\qbezier(180,80)(180,75)(170,75)
\qbezier(170,75)(160,75)(160,80)
\put(220,80){\vector(-1,0){10}}
\put(190,80){\vector(-1,0){10}}
\put(160,81){\vector(-1,1){12}}
\put(148,93){\vector(1,-2){16}}
\put(148,93){\circle{2}}
\put(164,61){\circle*{2}}
\put(220,80){\circle*{2}}
\put(210,80){\circle*{2}}
\put(190,80){\circle*{2}}
\put(180,80){\circle*{2}}
\put(160,81){\circle*{2}}
\put(139,93){\mbox{\footnotesize $\infty$}}
\put(165,63){\mbox{\footnotesize $1$}}
\put(204,78){\mbox{\footnotesize $3$}}
\put(222,78){\mbox{\footnotesize $2$}}
\put(192,78){\mbox{\footnotesize $4$}}
\put(174,78){\mbox{\footnotesize $5$}}
\put(162,78){\mbox{\footnotesize $6$}}
\qbezier(0,30)(0,50)(40,50)
\qbezier(80,30)(80,50)(40,50)
\qbezier(0,30)(0,10)(40,10)
\qbezier(80,30)(80,10)(40,10)
\qbezier(10,30)(10,35)(20,35)
\qbezier(30,30)(30,35)(20,35)
\qbezier(10,30)(10,25)(20,25)
\qbezier(30,30)(30,25)(20,25)
\qbezier(50,30)(50,35)(60,35)
\qbezier(70,30)(70,35)(60,35)
\qbezier(50,30)(50,25)(60,25)
\qbezier(70,30)(70,25)(60,25)
\put(40,50){\vector(-1,-1){16}}
\put(15,25){\vector(-1,-1){8}}
\put(56,34){\vector(-1,1){16}}
\put(73,17){\vector(-1,1){8}}
\put(40,50){\circle{2}}
\put(25,35){\circle*{2}}
\put(15,25){\circle*{2}}
\put(7,17){\circle*{2}}
\put(73,17){\circle*{2}}
\put(65,25){\circle*{2}}
\put(55,35){\circle*{2}}
\put(36,52){\mbox{\footnotesize $\infty$}}
\put(59,19){\mbox{\footnotesize $5$}}
\put(72,19){\mbox{\footnotesize $4$}}
\put(55,36){\mbox{\footnotesize $6$}}
\put(5,20){\mbox{\footnotesize $3$}}
\put(16,20){\mbox{\footnotesize $2$}}
\put(20,37){\mbox{\footnotesize $1$}}
\qbezier(140,30)(140,50)(180,50)
\qbezier(220,30)(220,50)(180,50)
\qbezier(220,30)(220,0)(200,0)
\qbezier(178,12)(190,0)(200,0)
\qbezier(140,30)(140,10)(160,10)
\qbezier(160,10)(180,10)(190,20)
\qbezier(190,20)(198,27)(207,27)
\qbezier(207,27)(213,27)(213,15)
\qbezier(213,15)(213,5)(200,5)
\qbezier(200,5)(190,5)(183,14)
\qbezier(174,18)(167,30)(160,30)
\qbezier(160,30)(153,30)(153,23)
\qbezier(153,23)(153,16)(160,16)
\qbezier(160,16)(177,16)(183,26)
\qbezier(183,26)(192,38)(180,38)
\qbezier(180,38)(175,38)(175,28)
\qbezier(175,28)(175,25)(177,21)
\put(180,50){\vector(-1,-1){20}}
\put(168,26){\vector(1,1){8}}
\put(187,33){\vector(2,-1){13}}
\put(203,27){\vector(-1,1){23}}
\put(180,50){\circle{2}}
\put(160,30){\circle*{2}}
\put(168,26){\circle*{2}}
\put(175,33){\circle*{2}}
\put(187,33){\circle*{2}}
\put(201,26){\circle*{2}}
\put(203,27){\circle*{2}}
\put(176,52){\mbox{\footnotesize $\infty$}}
\put(199,19){\mbox{\footnotesize $5$}}
\put(185,26){\mbox{\footnotesize $4$}}
\put(204,28){\mbox{\footnotesize $6$}}
\put(177,31){\mbox{\footnotesize $3$}}
\put(163,20){\mbox{\footnotesize $2$}}
\put(156,31){\mbox{\footnotesize $1$}}
\put(90,52){$h^\R_{2;3^1,\emptyset,\emptyset}=4$}
\end{picture}
%\centerline{\ig[scale=.15]{IMG_6.jpg}}
\caption{All the $4$ abstract diagrams of genus~$2$ with a single boundary pole of degree~$3$.
Orientation of the boundary arcs is not shown}\label{fd}
\end{figure}

Note that the orientation requirement implies that each internal pole, considered as a vertex
of the graph, has
an even degree (indeed, edges entering and leaving an internal pole must alternate).
Half of this degree is called the {\em order\/} of the internal pole.
(In other words, the order of an internal pole is its indegree, or, equivalently,
outdegree). The {\em order\/} of a boundary pole is one less than its degree.
Figure~\ref{fd} shows examples of nonisomorphic abstract diagram on a connected genus~$2$ surface. Each of these diagrams has $6$ critical points, one boundary pole of order~$3$ and no internal poles. Three of these diagrams are homeomorphic to a sphere with three discs removed while the fourth one is a torus without one disk.

\begin{remark}
If we erase in a diagram all internal edges that connect a critical point to a pole,
preserving only those connecting two critical points (as well as the orientations
and the numbering of all the critical values), then we can reconstruct the
original diagram in a unique way: inside each disk, connect the only pole on the
boundary by an internal edge with each critical point on the boundary
possessing no internal edges and orient the new edges properly. Sometimes it is more convenient to use
such reduced diagrams instead of complete ones.
\end{remark}

\subsection{The diagram of a framed real meromorphic function}

The diagram of a framed real meromorphic function $f:C\to\CP^1$ is, essentially, the graph embedded in the
domain~$C^f$ of the function~$f$ and formed by the preimage of the real line in the
target complex projective line. This graph is endowed with certain additional data
making it into a diagram.

\begin{definition}
Let $f:(C,\tau)\to(\CP^1,\sigma)$ be a framed simple real meromorphic function; here $\sigma:\CP^1\to\CP^1$
is complex conjugation, $\tau:C\to C$ is an anti-holomorphic involution, and
$f\circ\tau=\sigma\circ f$. Let~$m$ be the number of finite critical values of~$f$.
Number the finite critical values of~$f$ (which are all real) by numbers $1,2,\dots,m$,
successively in the increasing order, starting from the smallest one.
The {\em diagram} $D(f)$ of~$f$ consists of the following data:
\begin{itemize}
\item the domain $C^f$ of~$f$ endowed with the complex orientation;
\item the oriented graph in $C^f$, which is the preimage $f^{-1}(\RP^1)$,
whose vertices are the critical points and the poles of~$f$,
with the orientation of the edges induced by the natural orientation of the projective line $\RP^1\subset\CP^1$;
\item the numbering of the critical points at which the function~$f$ has a finite critical
value which associates to each critical point the number of the corresponding critical value.
\end{itemize}
\end{definition}
Clearly, the diagram of a framed simple real meromorphic function is an abstract diagram.
The orders of the poles in the diagram coincide with the orders of the poles of the function.

A theorem from~\cite{NSV02} establishes a one-to-one correspondence between
diagrams of framed real meromorphic functions and abstract diagrams. It follows that the enumeration of real meromorphic functions is reduced to the enumeration of their diagrams. For example, the four diagrams of Fig.~\ref{fd} contribute to the number $h^\R_{6;3^1,\emptyset,\emptyset}=4$ (where we assume that the orientation of these diagrams are induced from the standard orientation of the plane). Similar diagrams taken with their opposite orientations provide the computation of $h^\R_{6;\emptyset,3^1,\emptyset}=4$.

Note that the bidegree of a framed real meromorphic function introduced in Sec.~\ref{sFSHN}
can be read from the diagram of the function: it is formed by the numbers of the faces in the diagram
the orientation of whose boundaries is consistent or inconsistent with that of the underlying surface,
respectively.

\subsection{Proof of Theorem~\ref{tcjdpm}}\label{sP}

The coefficient $h^{\R\circ}_{m;(\kappa^+,\kappa^-,\lambda)}$ of the generating function $H^{\R\circ}$
is nothing but the number of diagrams with $m$ critical points
numbered from~$1$ to~$m$ and the partition~$\kappa^+$ (respectively,~$\kappa^-$) of orders of
positive (respectively, negative) real poles, and the partition~$\lambda$ of orders of internal poles.
Denote the set of such diagrams (both connected and disconnected ones) by $\cD^\circ_{m;(\kappa^+,\kappa^-,\lambda)}$.
Each diagram is counted with the coefficient inverse to the order of the automorphism
group of the diagram.

The proof of the theorem consists in establishing a one-to-one correspondence between the sets of diagrams with $m$ critical points and the set of diagrams with $m-1$ critical points and some additional markings. Introducing markings acts on generating functions as differential operators. The correspondence is based on considering local transformations of diagrams
that arise in the process of eliminating the critical point corresponding to the maximal critical
value. The differential operators corresponding to these local transformations are
the summands of the differential operator~$W^+$. For the differential operator~$W^-$,
the proof would be similar, but with the critical point corresponding to the maximal
critical value replaced by that for the minimal one. We prefer working with~$W^+$ because
eliminating the critical point with the maximal critical value allows one to preserve
the numbering of all the other critical points.

Take a diagram in $\cD^\circ_{m;(\kappa^+,\kappa^-,\lambda)}$.
This diagram contributes to the monomial $p^+_{\kappa^+}p^-_{\kappa^-}q_{\lambda}\frac{u^m}{m!}$
in the generating function~$H^{\R\circ}$.
There are two connected diagrams without critical points. Each of these diagrams is
a disc with a single pole, which can be either positive or negative.
These two diagrams correspond to the two framings of the linear function $z\mapsto z$.
Together with the (disconnected) diagram of Euler characteristic~$4$
added by our agreement in Sec.~\ref{sFSHN}, they provide the initial condition: for $u=0$, we have
$$
H^\R=\hbar^2(p^+_1+p^-_1)+\hbar^4q_1,\qquad H^{\R\circ}=e^{\hbar^2(p^+_1+p^-_1)+\hbar^4q_1}.
$$

Now suppose that~$m>0$, so that there is at least one critical point.
Consider the critical point number~$m$ of the diagram,
that is, the one corresponding to the largest finite critical value.
The internal edge issuing from this critical point can connect it
with either another critical point, or a pole. Below, we consider these two cases separately.

{\bf I. Critical point number~$m$ is connected by the internal edge with another critical
point.}

In this case there are at least~$2$ critical points, and
the internal edge is oriented to the critical point number~$m$, since any other critical point
has a smaller number. The two arcs of~$C^\tau$ issuing from the critical point number~$m$
are then outgoing,
and the other end of each of these arcs is a boundary pole. Once again, there are two cases:
the two poles on the ends of the arcs are distinct, or the two arcs have coinciding ends.

In the first case we replace the diagram by the one obtained from it by erasing the critical point number~$m$ and
by colliding the two boundary poles into a single one, placed instead of the critical point number~$m$.
The internal edge entering the critical point number~$m$ now enters the new pole instead.
The order of the new pole is equal to the sum of the orders of the two old poles.
One of the two old poles (the one to the right of the critical point number~$m$)
is necessarily positive, while the other one can be either positive or negative.
The new pole is positive or negative depending on whether the old pole to the left
is positive or negative. The new diagram contains $m-1$ critical points.
The orientation of all the edges is preserved.

\begin{figure}
\setlength{\unitlength}{0.5mm}
\thicklines
\begin{picture}(400,70)(0,0)
\qbezier(10,60)(0,40)(10,20)
\put(38,40){\vector(-1,0){33}}
\put(37,37){\vector(0,1){3}}
\put(37,43){\vector(0,-1){3}}
\qbezier(35,55)(40,40)(35,25)
\put(-5,38){\mbox{\footnotesize $m$}}
\put(40,38){\mbox{\footnotesize $m'$}}
\put(-4,53){\mbox{\footnotesize $\infty'$}}
\put(-7,23){\mbox{\footnotesize $\infty''$}}
\put(5,40){\circle*{2}}
\put(37,40){\circle*{2}}
\put(8,55){\circle*{2}}
\put(8,25){\circle*{2}}
\put(5,47){\vector(1,4){2}}
\put(5,33){\vector(1,-4){2}}
\put(8,55){\vector(2,-1){11}}
\put(19,61){\vector(-2,-1){11}}
\put(8,25){\vector(1,0){11}}
\put(55,40){\vector(1,0){12}}
\put(55,37){\line(0,1){6}}
\qbezier(90,60)(80,40)(90,20)
\put(118,40){\vector(-1,0){33}}
\qbezier(115,55)(120,40)(115,25)
\put(75,38){\mbox{\footnotesize $\infty$}}
\put(120,38){\mbox{\footnotesize $m'$}}
\put(117,37){\vector(0,1){3}}
\put(117,43){\vector(0,-1){3}}
\put(85,40){\circle*{2}}
\put(117,40){\circle*{2}}
\put(85,40){\vector(2,1){11}}
\put(91,49){\vector(-2,-3){6}}
\put(85,40){\vector(3,-2){9}}
\put(57,5){$(a)$}
\put(140,20){\mbox{\footnotesize $m$}}
\put(170,20){\mbox{\footnotesize $\infty$}}
\put(207,20){\mbox{\footnotesize $m'$}}
\put(150,25){\circle*{2}}
\put(205,25){\circle*{2}}
\put(85,40){\circle*{2}}
\put(170,25){\vector(0,1){12}}
\put(178,41){\vector(-1,-2){8}}
\qbezier(150,25)(160,10)(170,25)
%\put(160,25){\oval(20,15)[b]}
\qbezier(150,25)(150,26)(154,30)
\qbezier(156,32)(158,34)(160,34)
\qbezier(163,34)(166,34)(166,32)
\put(165,31){\vector(1,-1){6}}
\put(166,21){\vector(1,1){4}}
%\put(170,25){\line(-1,0){3}}
\qbezier(170,25)(177,40)(184,25)
\qbezier(150,25)(183,100)(205,25)
\qbezier(199,22)(209,25)(208,28)
\qbezier(192,20)(180,24)(186,27)
\qbezier(188,28)(190,29)(191,29)
%\put(150,25){\line(0,1){4}}
%\put(150,25){\line(1,0){4}}
%\put(150,25){\line(0,-1){4}}
%\put(154,21){\vector(-1,1){4}}
%\put(154,29){\vector(-1,-1){4}}
\put(154,33){\vector(-1,-2){4}}
\put(203,23){\vector(2,1){2}}
\put(207,27){\vector(-1,-1){2}}
\put(215,40){\vector(1,0){12}}
\put(215,37){\line(0,1){6}}
\put(215,5){$(b)$}
\qbezier(235,25)(255,100)(275,25)
\qbezier(235,25)(240,18)(250,18)
\qbezier(235,25)(235,26)(239,30)
\qbezier(241,32)(243,34)(245,34)
\qbezier(269,22)(279,25)(278,28)
\put(277,20){\mbox{\footnotesize $m'$}}
\put(240,57){\mbox{\footnotesize $\infty$}}
\put(275,25){\circle*{2}}
\put(250,59){\circle*{2}}
\put(250,59){\vector(1,-1){8}}
\put(244,47){\vector(1,2){6}}
%\put(250,59){\line(0,1){4}}
\put(253,62){\vector(-2,-1){4}}
\put(273,23){\vector(2,1){2}}
\put(277,27){\vector(-1,-1){2}}
\end{picture}
%\centerline{\ig[scale=.5]{IMG_7.jpg}}
\caption{The internal edge connects the critical point corresponding to the maximal critical value
with another critical point. Pictures show the diagram transformation in the cases when (a) the two neighboring
arcs connect the critical point to two distinct poles; (b) the two neighboring
arcs connect the critical point to one and the same pole}\label{fj}
\end{figure}
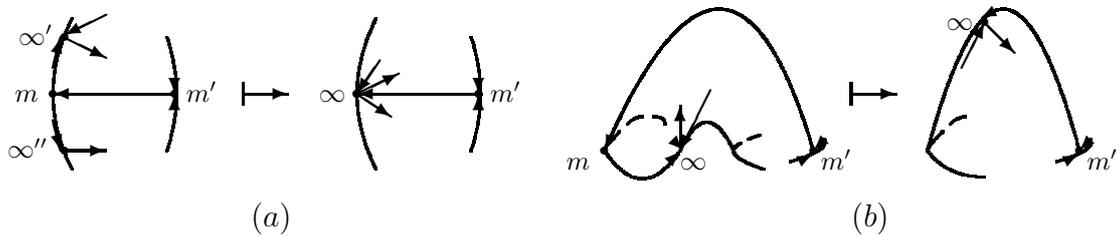

This local transformation of the diagram is shown in Fig.~\ref{fj} (a).
It causes the term
$$
\hbar p_{i+j}^{{\bar i}} \frac{\p^2}{\p p_i^{{\bar i}}\p p_j^+}
$$
in the differential operator~$W^+$. Indeed, a boundary pole of order~$i+j$ arises instead of two boundary poles,
of orders~$i$ and~$j$, respectively, and the sign of this pole coincides with that of the
pole to the left. These terms of the operator~$W^+$ form its ``join'' part, since they
join two cycles in the monodromy permutation over infinity, of lengths~$i$ and~$j$,
by a single cycle, of length~$i+j$.

If the two arcs issuing from the critical point number~$m$ end at the same boundary pole, then this pole is
necessarily of an even order. The critical point number~$m$ and the pole are the only two graph
vertices on the boundary component of the diagram containing them. We contract this
boundary component to a point, thus decreasing the genus of the diagram by~$1$,
make this point into an internal pole of order equal to half the order of the old boundary pole
and erase the critical point number~$m$, while preserving all the orientations of the internal edges.
As a result of this local transformation we obtain a diagram with~$m-1$ critical
points.

This local transformation of the diagram is shown in Fig.~\ref{fj} (b).
It causes the term
$$
\hbar q_i\frac\p{\p p^+_{2i}}
$$
in the differential operator~$W^+$. Indeed, the eliminated boundary pole has an even order~$2i$,
is positive, is replaced by a complex pole of order~$i$, and the genus of the diagram
decreases by~$1$.

{\bf II. Critical point number~$m$ is connected by the internal edge with a pole.}

This case has three subcases:
\begin{itemize}
\item $m=1$, so that there is a unique critical point;
\item $m>1$ and the critical point number~$m$ is connected by the internal edge with a boundary pole;
\item $m>1$ and the critical point number~$m$ is connected by the internal edge with an internal pole.
\end{itemize}

Let us consider these subcases one by one.

\begin{figure}
\setlength{\unitlength}{0.4mm}
\thicklines
\begin{picture}(400,70)(0,0)
\put(10,38){\circle{30}}
\put(-5,38){\vector(1,0){30}}
\put(-5,34){\vector(0,1){4}}
\put(-5,42){\vector(0,-1){4}}
\put(25,38){\circle{2}}
\put(-5,38){\circle*{2}}
\put(36,38){\vector(1,0){12}}
\put(36,35){\line(0,1){6}}
\put(65,56){\circle{30}}
\put(80,56){\circle{2}}
\put(65,20){\circle{30}}
\put(80,20){\circle{2}}
\qbezier(50,56)(65,45)(80,56)
\qbezier(50,56)(51,59)(54,59)
\qbezier(57,61)(58,62)(61,62)
\qbezier(64,62)(66,63)(68,62)
\qbezier(72,62)(75,62)(76,61)
\qbezier(78,59)(79,59)(80,56)
\qbezier(50,20)(65,9)(80,20)
\qbezier(50,20)(51,23)(54,23)
\qbezier(57,25)(58,26)(61,26)
\qbezier(64,26)(66,27)(68,26)
\qbezier(72,26)(75,26)(76,25)
\qbezier(78,23)(79,23)(80,20)
\put(26,35){\mbox{\footnotesize $\infty$}}
\put(-10,35){\mbox{\footnotesize $1$}}
\put(82,53){\mbox{\footnotesize $\infty'$}}
\put(82,17){\mbox{\footnotesize $\infty''$}}
\put(30,5){$(a)$}
\qbezier(110,63)(105,38)(110,13)
\qbezier(141,53)(146,38)(141,23)
\put(108,38){\vector(1,0){36}}
\put(107,38){\circle*{2}}
\put(144,38){\circle{2}}
\put(107,34){\vector(0,1){4}}
\put(107,42){\vector(0,-1){4}}
\put(99,36){\mbox{\footnotesize $m$}}
\put(146,36){\mbox{\footnotesize $\infty$}}
\put(108,52){\circle*{2}}
\put(108,24){\circle*{2}}
\put(98,50){\mbox{\footnotesize $m'$}}
\put(98,22){\mbox{\footnotesize $m''$}}
\put(124,60){\vector(-2,-1){16}}
\put(124,16){\vector(-2,1){16}}
\put(144,38){\vector(-1,1){14}}
\put(144,38){\vector(-1,-1){14}}
\put(155,38){\vector(1,0){12}}
\put(155,35){\line(0,1){6}}
\qbezier(180,20)(160,30)(200,30)
\qbezier(220,20)(240,30)(200,30)
\qbezier(180,56)(160,46)(200,46)
\qbezier(220,56)(240,46)(200,46)
\put(174,26){\circle*{2}}
\put(174,50){\circle*{2}}
\put(164,24){\mbox{\footnotesize $m''$}}
\put(164,48){\mbox{\footnotesize $m'$}}
\put(226,26){\circle{2}}
\put(226,50){\circle{2}}
\put(229,24){\mbox{\footnotesize $\infty''$}}
\put(229,48){\mbox{\footnotesize $\infty'$}}
\put(189,26){\vector(-1,0){15}}
\put(189,50){\vector(-1,0){15}}
\put(226,26){\vector(-1,0){15}}
\put(226,50){\vector(-1,0){15}}
\put(157,5){$(b)$}
\qbezier(250,40)(260,20)(275,15)
\put(260,25){\circle*{2}}
\put(250,23){\mbox{\footnotesize $m$}}
\qbezier(260,25)(270,45)(280,60)
\put(280,60){\circle{2}}
\put(268,57){\mbox{\footnotesize $\infty$}}
\put(276,54){\vector(1,2){4}}
\put(280,60){\vector(-1,1){12}}
\put(292,72){\vector(-1,-1){12}}
\put(280,60){\vector(1,0){16}}
\put(292,48){\vector(-1,1){12}}
\put(280,60){\vector(0,-1){16}}
\put(257,28){\vector(1,-1){4}}
\put(263,22){\vector(-1,1){4}}
\qbezier(320,40)(330,20)(345,15)
\put(294,38){\vector(1,0){12}}
\put(294,35){\line(0,1){6}}
\put(330,25){\circle{2}}
\put(316,23){\mbox{\footnotesize $\infty$}}
\put(330,25){\vector(1,1){12}}
\put(330,25){\vector(0,1){16}}
\put(346,25){\vector(-1,0){16}}
\put(330,25){\vector(2,-1){16}}
\put(338,41){\vector(-1,-2){8}}
\put(327,28){\vector(1,-1){4}}
\put(333,22){\vector(-1,1){4}}
\put(295,5){$(c)$}
\end{picture}
%\centerline{\ig[scale=.5]{IMG_8.jpg}}
\caption{The internal edge connects the critical point corresponding to the maximal critical value
with a pole. Pictures show the diagram transformation in the cases when (a)
$m=1$; (b) $m>1$ and the internal edge connects the critical point to a boundary pole;
(c) $m>1$ and the internal edge connects the critical point to an internal pole}\label{fc}
\end{figure}
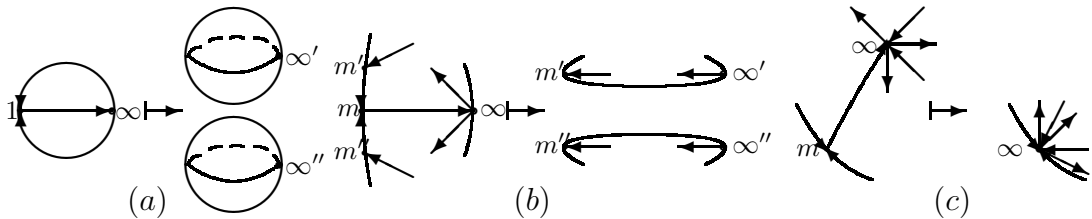

If $m=1$, then there exists a single connected diagram. The genus of the diagram is~$0$,
the function is $z\mapsto z^2$,
see Fig~\ref{fc} (a). The corresponding local
transformation consists in replacing this connected component
of a diagram with a pair of connected components, each being~$\CP^1$
with a single boundary pole. This operation decreases the genus of the diagram by~$1$ and leads to the summand
$$
p_1^{{-}}p_1^+\frac\p{\p p_{2}^{{-}}}
$$
in the operator~$W^+$.

If the critical point number~$m>1$ is connected by an internal edge with a boundary pole,
then both neighboring vertices of this critical point along the boundary arcs of the diagram are
critical points (which may well coincide): otherwise there would be a face of the diagram whose
boundary passes through poles at least twice.
The local transformation erasing the critical point number~$m$ splits the boundary pole into two boundary poles,
whose orders sum to the order of the original pole, see Fig.~\ref{fc} (b).
This operation decreases the genus of the diagram by~$1$ and leads to the summand
$$
p_i^{{\bar i}}p_j^+\frac\p{\p p_{i+j}^{{\bar i}}}
$$
in the operator~$W^+$. This is the ``cut'' part of the cut-and-join operator.

Indeed, the sign of the first new pole coincides with that of the original one, while
the second pole is necessarily positive, whatever is the sign of the original pole.

If the internal edge issuing from the critical point number~$m$ connects it with an internal
pole of order~$i$, then the local transformation consists in contracting this edge,
see Fig.~\ref{fc} (c). The internal pole is replaced by a boundary pole of order~$2i$
and the genus of the diagram is preserved, which leads to the summand
$$
p_{2i}^+\frac\p{\p q_i}
$$
in the operator~$W^+$. This completes the proof of Theorem~\ref{tcjdpm}.

Note that each local transformation of a diagram used in the proof
establishes a one-to-one correspondence between the faces in the original
and the resulting diagrams, as well as the orientation of the boundary of each face.
This is just another way to show that the operator~$W^+$ respects both the
grading and the bigrading in the space of polynomials.

\section{A topological field theory description}\label{sTFT}

In this section we describe yet another combinatorial reformulation of the diagrammatic
model used in Section~\ref{sDD}. For notation simplicity, we set $\hbar=1$ throughout this section. The necessary rescaling of the variables corresponding to insertion of $\hbar$ in the general case, follows easily from Riemann-Hurwitz formula, see Remark~\ref{remh}. The description of this section is similar to the well-known
interpretation of ordinary complex Hurwitz numbers as correlators
of a toy topological field theory (see~\cite{D95,AN06,AN08-1,AN08-2}),
which can be formulated as follows. Let~$V_n$ denote the center of the
group algebra of  the symmetric group~$S_n$, $V_n=Z\C[S_n]$.
This vector space is freely spanned by the conjugacy classes~$v_\mu$ of permutations in~$S_n$,
these conjugacy classes being numbered by partitions $\mu\vdash n$.
Multiplication inherited by~$V_n$ from the group algebra $\C[S_n]$ makes it
into a Frobenius algebra, with the scalar product $(a,b)$ defined as the
coefficient of the conjugacy class $[{\rm id}]$ of the identity permutation in the product $ab$, so that $(ab,c)=(a,bc)$
for any triple $a,b,c\in V_n$.

The topological field theory in question assigns, for a given orientable surface $\Sigma$ with~$m$ punctures,
the vector space~$V_n$ to each of the punctures. If we pick an $m$-tuple $(\mu_1,\dots,\mu_m)$
of partitions of~$n$, then the coefficient of~$[{\rm id}]$ in the product $v_{\mu_1}\dots v_{\mu_m}$
enumerates degree~$n$ ramified coverings of~$\Sigma$ that are unramified out of~$m$
marked points and have prescribed ramification types $\mu_1,\mu_2,\dots,\mu_m$
over the marked points. If the surface~$\Sigma$ is the sphere, and the ramification type
over one of the marked points is~$\mu$, while the ramification type over all the other
marked points is~$1^{n-2}2^1$, then these ramified coverings are enumerated by simple Hurwitz numbers.

Framed real Hurwitz numbers admit a similar description, which we present below.

\subsection{States and transitions}

In this section we refine the construction of an algebra of pairs of involutions from~\cite{AN06}, Sec.~5.
The refinement adjusts this construction to the case of simple purely real framed Hurwitz numbers.
Pick a finite set $N$, $n=|N|$ being the number of elements in~$N$,
and a representation of~$N$ as a disjoint union $N=N^+\bigsqcup N^-$
 of two subsets $N^+, N^-$ consisting of~$n^+$ and~$n^-$ elements, respectively, $n^++n^-=n$.

\begin{definition} A \emph{state} is a partition of~$N$ into a disjoint union of one and two-element subsets such that each two-element subset contains one element from $N^+$ and one from~$N^-$. A \emph{transition} is an ordered pair of states;
the first state in a transition is said to be {\em initial}, the second one {\em final}.
\end{definition}

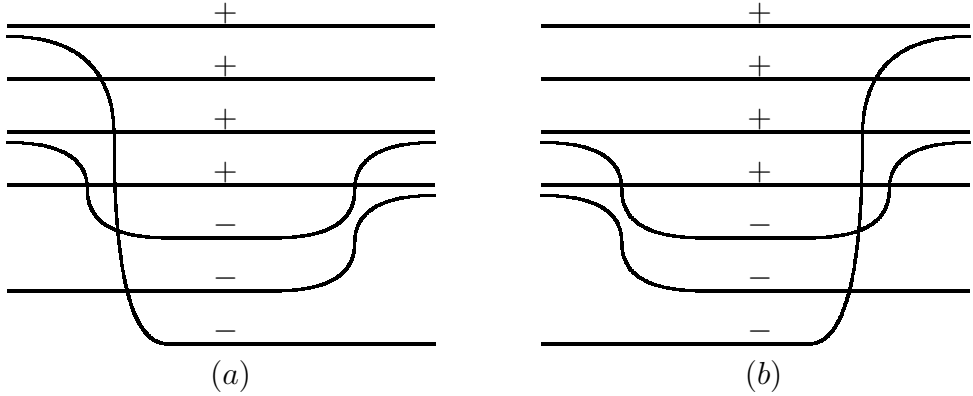
\begin{figure}
\thicklines
\begin{picture}(400,120)(100,-20)
\multiput(90,60)(0,20){4}{\line(1,0){160}}
\multiput(167,62)(0,20){4}{$+$}
\multiput(150,0)(0,20){3}{\line(1,0){40}}
\multiput(167,2)(0,20){3}{$-$}
\qbezier(150,40)(120,40)(120,58)
\qbezier(120,58)(120,76)(90,76)
\put(90,20){\line(1,0){60}}
\put(190,0){\line(1,0){60}}
\qbezier(150,0)(130,0)(130,78)
\qbezier(130,78)(130,116)(90,116)
\qbezier(190,40)(220,40)(220,58)
\qbezier(220,58)(220,76)(250,76)
\qbezier(190,20)(220,20)(220,38)
\qbezier(220,38)(220,56)(250,56)
\put(166,-15){$(a)$}
\multiput(290,60)(0,20){4}{\line(1,0){160}}
\multiput(367,62)(0,20){4}{$+$}
\multiput(350,0)(0,20){3}{\line(1,0){40}}
\multiput(367,2)(0,20){3}{$-$}
\qbezier(390,40)(420,40)(420,58)
\qbezier(420,58)(420,76)(450,76)
\put(390,20){\line(1,0){60}}
\put(290,0){\line(1,0){60}}
\qbezier(390,0)(410,0)(410,78)
\qbezier(410,78)(410,116)(450,116)
\qbezier(350,40)(320,40)(320,58)
\qbezier(320,58)(320,76)(290,76)
\qbezier(350,20)(320,20)(320,38)
\qbezier(320,38)(320,56)(290,56)
\put(366,-15){$(b)$}
\end{picture}
\caption{(a) A transition on a set~$N=N^+\sqcup N^-$, $|N^+|=n^+=4$, $|N^-|=n^-=3$. The initial state is
represented by the involution on the left part of the picture, while
the final state by that on the right. An orbit of length two in each involution is depicted by
an arc issuing from the negative element of the orbit
and approaching the corresponding positive element from below. The type of
the transition is $(1^12^1,2^1,1^1)$. (b) Its inverse, of type $(2^1,1^12^1,1^1)$.
}\label{ftrans}
\end{figure}

Thus, for any two states there is a unique transition from the first of them to the
second one. An example of a transition is depicted in Fig.~\ref{ftrans}~(a).
The points of $N$ correspond to the horizontal arcs. They are marked with signs $+$ or $-$ depending on whether they belong to $N^+$ or $N^-$, respectively. The initial and the final states of the transition correspond to the left and right parts of the picture, respectively: the points forming a single pair of a state are depicted as two approaching arcs.

\begin{definition} The \emph{ type } of a transition is its orbit under the action of the group $S(n^+)\times S(n^-)$
acting on the set of transitions by separately permuting the elements in $N^+$ and $N^-$.
\end{definition}

We label transition types by triples of partitions $\mu=(\kappa^+,\kappa^-,\lambda)=((k_1^+,k_2^+,\dots),(k_1^-,k_2^-,\dots),(\ell_1,\ell_2,\dots))$ in the following way. Pairs of the initial and the final states determine chains of elements of the form: an element $x_1$, its couple $x_2$ in the initial state, the couple~$x_3$ of $x_2$ in the final state, the couple~$x_4$ of~$x_3$ in the initial state, etc. Each such chain is either cyclic or it starts and ends in the elements having no couples in one of the two states. We order elements of a chain in such a way that a positive point of a pair follows a negative one if it is a pair in the initial state and a positive point precedes a negative one if they are points of one pair in the final state. According to that we define the \emph{sign of a chain} which is not a cycle as the sign of its first element. More explicitly, if this chain has an odd number of elements, then its sign is the sign of either of its ends. If the chain has an even number of elements, its sign is positive or negative if both its ends are single elements of the initial or the final state, respectively.
%The signs of the chains contained in the transition of Fig.~\ref{ftrans} are shown to the right of the picture.

The \emph{transition type} is equal to
$$\mu=((k_1^+,k_2^+,\dots),(k_1^-,k_2^-,\dots),(\ell_1,\ell_2,\dots))$$
 if the whole set $N$ splits
 into a disjoint union of positive chains of lengths $k_1^+,k_2^+,\dots$, negative chains of lengths $k_1^-,k_2^-,\dots$, and cyclic chains of lengths $2\ell_1,2\ell_2,\dots$.

Remark that the numbers $n^+$ and $n^-$ can be recovered from the transition type uniquely. Namely,
\begin{equation}\label{npm}
\begin{aligned}
n^+(\kappa^+,\kappa^-,\lambda)&=\sum_{k\in \kappa^+}\lceil k/2\rceil+\sum_{k\in \kappa^-}\lfloor k/2\rfloor+\sum_{\ell\in \lambda}\ell,\\
n^-(\kappa^+,\kappa^-,\lambda)&=\sum_{k\in \kappa^+}\lfloor k/2\rfloor+\sum_{k\in \kappa^-}\lceil k/2\rceil+\sum_{\ell\in \lambda}\ell,\\
\end{aligned}
\end{equation}
where $\lfloor\cdot\rfloor$ and $\lceil\cdot\rceil$ denote the floor and the ceiling functions, respectively.

The \emph{inverse transition} is obtained from the given one by exchanging its initial and final states.
Fig.~\ref{ftrans}~(b) shows the inverse transition to the one in Fig.~\ref{ftrans}~(a).
 If the original transition is of type $\mu=(\kappa^+,\kappa^-,\lambda)$, then the type of its inverse is obtained from $\mu$ by exchanging parts of even lengths in the partitions $\kappa^+$ and $\kappa^-$.

A transition is said to be \emph{trivial} if the initial and the final states are equal. The type of a trivial transition has the form $(1^i,1^j,1^{k})$ for some $i,j,k$ with $i+k=n^+$, $j+k=n^-$.

A transition is called a \emph{transposition } if the initial and the final states differ by two elements which are joined in a pair in one of the states and are individual elements of the other state. The type of a transposition is either $(1^{i-2}2^1,1^j,1^{k})$ or $(1^i,1^{j-2}2^1,1^{k})$. These transposition
types are said to be \emph{positive} and \emph{negative}, respectively.

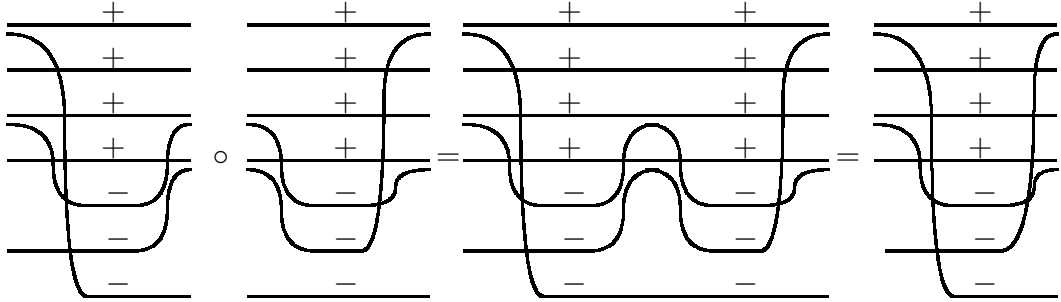
\begin{figure}
\setlength{\unitlength}{0.3mm}
\thicklines
\begin{picture}(400,120)(20,0)
\multiput(40,60)(0,20){4}{\line(1,0){80}}
\multiput(81,62)(0,20){4}{$+$}
\multiput(75,0)(0,20){3}{\line(1,0){20}}
\multiput(83,2)(0,20){3}{$-$}
\qbezier(75,40)(60,40)(60,58)
\qbezier(60,58)(60,76)(40,76)
\put(40,20){\line(1,0){40}}
\put(90,0){\line(1,0){30}}
\qbezier(75,0)(65,0)(65,78)
\qbezier(65,78)(65,116)(40,116)
\qbezier(95,40)(110,40)(110,58)
\qbezier(110,58)(110,76)(120,76)
\qbezier(95,20)(110,20)(110,38)
\qbezier(110,38)(110,56)(120,56)
\put(130,58){$\circ$}
\multiput(145,60)(0,20){4}{\line(1,0){80}}
\multiput(183,62)(0,20){4}{$+$}
\multiput(175,0)(0,20){3}{\line(1,0){20}}
\multiput(183,2)(0,20){3}{$-$}
\qbezier(195,40)(210,40)(210,48)
\qbezier(210,48)(210,56)(225,56)
\put(195,0){\line(1,0){30}}
\put(145,0){\line(1,0){30}}
\qbezier(195,20)(205,20)(205,88)
\qbezier(205,88)(205,116)(225,116)
\qbezier(175,40)(160,40)(160,58)
\qbezier(160,58)(160,76)(145,76)
\qbezier(175,20)(160,20)(160,38)
\qbezier(160,38)(160,56)(145,56)
\put(228,58){$=$}
\multiput(240,60)(0,20){4}{\line(1,0){80}}
\multiput(281,62)(0,20){4}{$+$}
\multiput(275,0)(0,20){3}{\line(1,0){20}}
\multiput(283,2)(0,20){3}{$-$}
\qbezier(275,40)(260,40)(260,58)
\qbezier(260,58)(260,76)(240,76)
\put(240,20){\line(1,0){40}}
\put(295,0){\line(1,0){30}}
\qbezier(275,0)(265,0)(265,78)
\qbezier(265,78)(265,116)(240,116)
\qbezier(295,40)(310,40)(310,58)
\qbezier(310,58)(310,76)(325,76)
\qbezier(295,20)(310,20)(310,38)
\qbezier(310,38)(310,56)(325,56)
\multiput(320,60)(0,20){4}{\line(1,0){80}}
\multiput(358,62)(0,20){4}{$+$}
\multiput(350,0)(0,20){3}{\line(1,0){20}}
\multiput(358,2)(0,20){3}{$-$}
\qbezier(370,40)(385,40)(385,48)
\qbezier(385,48)(385,56)(400,56)
\put(370,0){\line(1,0){30}}
\put(320,0){\line(1,0){30}}
\qbezier(370,20)(380,20)(380,88)
\qbezier(380,88)(380,116)(400,116)
\qbezier(350,40)(335,40)(335,58)
\qbezier(335,58)(335,76)(320,76)
\qbezier(350,20)(335,20)(335,38)
\qbezier(335,38)(335,56)(320,56)
\put(403,58){$=$}
\multiput(420,60)(0,20){4}{\line(1,0){80}}
\multiput(461,62)(0,20){4}{$+$}
\multiput(455,0)(0,20){3}{\line(1,0){20}}
\multiput(463,2)(0,20){3}{$-$}
\qbezier(455,40)(440,40)(440,58)
\qbezier(440,58)(440,76)(420,76)
\put(425,20){\line(1,0){30}}
\put(475,0){\line(1,0){25}}
\qbezier(455,0)(445,0)(445,78)
\qbezier(445,78)(445,116)(420,116)
\qbezier(475,40)(490,40)(490,48)
\qbezier(490,48)(490,56)(500,56)
\qbezier(475,20)(490,20)(490,88)
\qbezier(490,88)(490,116)(500,116)
\end{picture}
\caption{A nonzero composition of two transitions
}\label{ftcomp}
\end{figure}

A composition of two transitions $(s_1,s'_1)$ and $(s_2,s'_2)$
exists if and only if the final state of the first transition is equal to the initial state of the second one,
$s'_1=s_2$; in this case the composition equals $(s_1,s'_2)$.
In particular, the composition of a transition and its inverse is a trivial transition.
An example of a nonzero composition of two transitions is shown in Fig.~\ref{ftcomp}.

Now let~$f:(C,\tau)\to(\CP^1,\sigma)$ be a framed real function.
The total preimage~$f^{-1}(\RP^1)$ of the real line splits the
domain~$C$ of~$f$ into open discs.
Take for the set~$N_f$ associated to this function the set of discs
that are the preimages of the upper hemisphere.
The number~$n_f=|N_f|$ of elements in~$N_f$ is the degree of~$f$.
Each of these discs is a subset in either~$C^f$, or its complement $C\setminus C^f$,
which determines the splitting of~$N_f$ into two subsets~$N_f^+$, $N_f^-$.
Complex conjugation~$\tau$ establishes a one-to-one correspondence
between the discs in~$N_f$ and the discs that are preimages of the lower hemisphere.

The finite critical values of~$f$, together with infinity, cut the real projective line~$\RP^1\subset\CP^1$
into~$m+1$ arcs. The arc between $k$~th and~$k+1$~th critical values determines
a state of the set~$N_f$: two discs in~$N_f$ belong to the same pair iff
the common boundary between the first of them and the $\tau$-conjugate
of the second one contains an arc connecting critical points number~$k$ and~$k+1$.
Hence, each critical value determines a transition between two consecutive states.
The transition corresponding to a simple critical value is a transposition.

\begin{proposition}\label{proptr} For any triple $\mu$ of partitions
the Hurwitz number $h^{\R\circ}_{m;\mu}$ enumerating not necessarily connected
framed simple real meromorphic functions is equal to the number of sequences of~$m$ transpositions
such that their composition is defined and has the given type $\mu$,
divided by the factor $n^+!n^-!$,
\begin{equation}\label{tauseq}
h^{\R\circ}_{m;\mu}=\frac{1}{n^+!n^-!}\left|\left\{(\tau_1,\dots,\tau_m)\Bigm| \tau_i\hbox{ is a transposition, }
    \tau_1\dots\tau_m\hbox{ is of type }\mu\right\}\right|.
\end{equation}
\end{proposition}

\emph{Proof.}
The set~$N$  of faces of any diagram in the sense of Sect.~\ref{sDD} is split into two subsets,
$N^+$ and~$N^-$ according to whether the orientation of the boundary of the face
coincides with the one induced by the orientation of the face.
It is sufficient to establish a bijection between all possible sequences of transitions
contributing to the right hand side of Eq.~\eqref{tauseq} and the set of diagrams
with numbered faces.

Let $\tau_1,\dots,\tau_m$ be a sequence of transitions. Let $s_0,\dots,s_m$ be the corresponding sequence of states so that $\tau_k$ is the transition from $s_{k-1}$ to $s_k$. Consider also the sequence of segments $I_0=[\infty,1]$, $I_1=[1,2]$, \dots, $I_{m-1}=[1,m]$, $I_m=[m,\infty]$ of the real projective line $\R P^1=\R\bigcup\{\infty\}$. The corresponding diagram is defined as the union of discs $D_i$ labeled by the indices $i\in N$. Each disk is isomorphic to the upper half-sphere $\{{\rm Im}z\ge0\}\subset\C P^1$. The discs~$D_i$ and~$D_j$ are glued along the segment~$I_k$ if the points~$i$ and~$j$ belong to one pair in the state~$s_k$. It is easy to verify all necessary properties of this correspondence.

To be more precise, the ramification type at infinity of the meromorphic function corresponding to the obtained diagram is equal to the type of the transition between $s_m$ and $s_0$ which is \emph{inverse} to the transition type $\mu$ of the composition $\tau_1\dots\tau_m$. But the Hurwitz numbers
$h^{\R\circ}_{m;(\kappa^+,\kappa^-,\lambda)}$ and $h^{\R\circ}_{m;(\kappa^-,\kappa^+,\lambda)}$ are equal, due to symmetry. The proposition follows.

\subsection{The Frobenius algebras of transitions}

One can reformulate the statement of Proposition~\ref{proptr} as follows. Define the algebra $T_{n^+,n^-}$ (say, over $\Q$) spanned by all possible transitions for a given set $N=N^+\bigsqcup N^-$, $|N^\pm|=n^\pm$,
with the product defined by the composition (if the composition of two transitions is not defined we set the product equal to zero).
Define $A_{n^+,n^-}=T_{n^+,n^-}^{S(n^+)\times S(n^-)}$ to be the $S(n^+)\times S(n^-)$-invariant part of $T_{n^+,n^-}$.
Denote by $C_\mu\in T_{n^+,n^-}$ the sum of all transitions of type~$\mu$.
Then $A_{n^+,n^-}$ is an associative commutative algebra with the basis $C_\mu$ for all triples of partitions $\mu=(\kappa^+,\kappa^-,\lambda)$ with $n^\pm(\mu)=n^\pm$. The unit of $A_{n^+,n^-}$ is $1\!\!1_{n^+,n^-}=\sum_{{i+k=n^+\atop j+k=n^-}}C_{1^i,1^j,1^k}$, the sum of all trivial transitions. Then the assertion of Proposition~\ref{proptr} can be reformulated as follows: the Hurwitz number $h^{\R\circ}_{m;\mu}$ is a suitably rescaled coefficient of $C_\mu$ in $C_2^m$ in the algebra $A_{n^+,n^-}$,
\begin{equation}\label{eC2m}
\frac{C_2^m\,1\!\!1_{n^+,n^-}}{n^+!n^-!}=\sum_{\mu,\;n^\pm(\mu)=n^\pm} h^{\R\circ}_{m;\mu}\frac{C_\mu}{|C_\mu|},
\end{equation}
where $C_2$ is the sum of all simple transitions, and $|C_\mu|$ is the number of distinct transitions of type $\mu$ in $T_{n^+,n^-}$.

We can identify now $A_{n^+,n^-}$ with the space of bidegree $(n^+,n^-)$-homogeneous polynomials in the variables
$p_k^+$, $p_k^-$, $q_k$, $k=1,2,\dots$ with the bidegree given by~Eq.~(\ref{ebd})
%$\bdeg~ p_{2k}^\pm=\bdeg~ q_k=(k,k),\quad\bdeg~ p_{2k+1}^+=(k+1,k),\quad\bdeg~ p_{2k+1}^-=(k,k+1)$
 by assigning
\begin{equation}\label{epC}
p_\mu\longleftrightarrow \frac{C_\mu}{|C_\mu|},
\end{equation}
where for $\mu=((k_1^+,\dots),(k_1^-,\dots),(\ell_1,\dots))$ we set $p_\mu=\prod p_{k_i^+}^+\prod p_{k_i^-}^-\prod q_{\ell_i}$.

\begin{theorem}\label{thC2m}
With the identification~{\rm \eqref{epC}}, the operator of multiplication on the left by $C_2$ in $A_{n^+,n^-}$ acts on polynomials in
the variables $p^\pm_k,q_k$ as the differential operator $W^+$ of Theorem~{\rm \ref{tcjdpm}}.
\end{theorem}

With the equality of this theorem, Eq.~\eqref{eC2m} becomes
$$(W^+)^m\sum_{i+k=n^+,\,j+k=n^-}\frac{(p_1^+)^i}{i!}\frac{(p_1^-)^j}{j!}\frac{q_1^k}{k!}=\sum_{\mu,\;n^\pm(\mu)=n^\pm} h^{\R\circ}_{m;\mu}p_\mu,$$
or, summing over all $n^+,n^-$, and $m$ with the factor $u^m/m!$ we obtain a different proof of the equality
$$H^{\R\circ}=e^{u W^+}e^{p_1^++p_1^-+q_1}$$
implying, in particular, the differential equation of Theorem~\ref{tcjdpm}
$$\frac{\p H^{\R\circ}}{\p u}=W^+ H^{\R\circ}.$$

\emph{Proof of Theorem}~\ref{thC2m}. Let us represent $C_2=C_2^{\rm cut}+C_2^{\rm join}$ where $C_2^{\rm cut}$ and $C_2^{\rm join}$ are the sums of positive and negative transpositions, respectively.

Then $C_2^{\rm join}$ acts on transitions by joining two individual elements of the initial state of a transition to a pair.
If these two elements belong to a single chain of length $2k$, then this chain is positive and the result of joining
is a cyclic chain of length~$2k$. This possibility corresponds to the summand $\sum q_{k}\frac{\p}{\p p_{2k}^+}$ of $W^+$.
If the positive and the negative joined elements belong to different chains of lengths~$j$ and $i$, respectively,
then these chains are joined to a single chain of length~$i+j$. Moreover, the sign convention implies
that the signs of the original two chains and the resulting one are~$+$, $\overline i$ and~$\overline i$, respectively.
This possibility corresponds to the summand $\sum p_{i+j}^{\overline i}\frac{\p^2}{\p p_i^{\overline i}p_{j}^+}$ of $W^+$.

The operator $C_2^{\rm cut}$ acts on transitions by cutting a pair of elements of the initial state of transition into two individual elements of the state.
If these two elements belong to a cyclic chain of length $2k$, then the result of cutting is a (positive) non-cyclic chain of length~$2k$.
This possibility corresponds to the summand $\sum k p_{2k}^+\frac{\p}{\p q_{k}}$ of $W^+$.
The factor $k$ reflects the fact that the original cyclic chain of length $2k$ contains exactly $k$ pairs of elements to which such cut can be applied.
Finally, if the elements of the pair belong to a non-cyclic chain, then the result of cutting are two chains. If we denote by~$i$ and~$j$ the lengths of these two chains containing the negative and the positive elements of the pair, respectively, then the signs of the initial and the two resulting chains are~$\overline i$, $\overline i$, and~$+$, respectively. This possibility corresponds to the summand
$\sum p_{i}^{\overline i}p_j^+\frac{\p}{\p p_{i+j}^{\overline i}}$ of $W^+$. This completes the proof of Theorem~\ref{thC2m}. The details are left to the reader.

\begin{remark}
Similarly to the operator~$W^+$, which corresponds to multiplication on the left by~$C_2$ in the algebra~$A_{n^+,n^-}$,
the operator~$W^-$ corresponds to multiplication by~$C_2$ on the right. This observation explains why the two
operators commute.
\end{remark}

Note that the bigrading of the variables imply the following product representation for the generating function
\begin{eqnarray*}
\cA(x,y)&=&\sum_{n^\pm=0}^\infty \dim~A_{n^+,n^-}x^{n^+}y^{n^-}\\
&=&1+(x+y)+(x^2+4xy+y^2)+(x^3+5x^2y+5xy^2+y^3)\\
&&\quad+(x^4+5x^3y+15x^2y^2+5xy^3+y^4)\\
&&\quad+(x^5+5x^4y+19x^3y^2+19x^2y^3+5xy^4+y^5)+\dots
\end{eqnarray*}
for the dimensions of the vector spaces~$A_{n^+,n^-}$:
$$
\cA(x,y)=\frac1{\prod_{k=1}^\infty(1-x^ky^k)^3(1-x^ky^{k-1})(1-x^{k-1}y^k)}.
$$

\medskip
Now we are going to introduce in each vector space~$A_{n^+,n^-}$ a nondegenerate scalar product
with respect to which the operators~$W^+$, $W^-$, and~$W$ are self-adjoint.
For a transition type represented by a triple of partitions $\mu=(\kappa^+,\kappa^-,\lambda)$ we denote by $\zeta(\mu)$ the cardinality of the stabilizer in $S(n^+)\times S(n^-)$ of an arbitrary transition of type $\mu$, so that $|C_\mu|=\frac{n^+!n^-!}{\zeta(\mu)}$. More explicitly,
$$
\zeta(\mu)=|\Aut(\kappa^+)|\;|\Aut(\kappa^-)|\;|\Aut(\lambda)|\;\prod_{j=1}^{\ell(\lambda)}\ell_j,
$$
where for a partition $\lambda$, denote by $|\Aut(\lambda)|$ the cardinality of its group of automorphisms, that is,
the product of factorials of the numbers of repeating parts.
The scalar product in~$A_{n^+,n^-}$ is defined by
\begin{equation}\label{esp}
(C_\mu,C_\nu)=\delta_{\mu,\nu}\;n^+!\;n^-!\;|C_\mu|.
\end{equation}
Together with the isomorphism~\eqref{epC} it defines the corresponding scalar product in the space of polynomials in
the variables $(p_k^\pm,q_k)$:
$$
(p_\mu,p_\nu)=\delta_{\mu,\nu}\;\zeta(\mu).
$$

\begin{proposition}\label{propsa}
The operators~$W^+,W^-$ on each of the vector spaces~$A_{n^+,n^-}$ are self-adjoint with respect to
scalar product~{\rm\eqref{esp}}.
\end{proposition}

In addition, the operators $W^+$ and $W^-$ commute. As an immediate corollary we obtain

\begin{corollary}
In each vector space~$A_{n^+,n^-}$, the operators~$W^+$, $W^-$, and $W$ admit a basis consisting of common
eigenvectors.
\end{corollary}

Properly chosen elements of these bases can be considered as the ``real'' analogues of the Schur polynomials.
For example, for $(n^+,n^-)=(1,1)$ the vector space $A_{1,1}$ is $4$-dimensional,
spanned by the monomials $p_2^+,p_2^-,p_1^+p_1^-,q_1$,
the eigenbasis is unique up to multiplication by constants and has the form
\begin{eqnarray*}
p_2^++p_2^-+p_1^+p_1^-+q_1,\qquad p_2^++p_2^--p_1^+p_1^-+q_1,&&\\
p_2^+-p_2^-+p_1^+p_1^--q_1,\qquad
p_2^+-p_2^--p_1^+p_1^--q_1.
\end{eqnarray*}
These polynomials govern the representation theory of the algebras~$A_{n^+,n^-}$, and their study is the subject
of further investigation.

\medskip
\emph{Proof of Proposition}~\ref{propsa}. Denote by $w^+_{\mu\nu}$ the matrix coefficients of the operator $W^+$: $W^+C_\nu=\sum_\mu w^+_{\mu\nu}C_\mu$. By definition, $w^+_{\mu\nu}\,|C_\mu|$ is equal to the number of triples $(a,b,\tau)$ of transitions such that $a$ and $b$ have transition types $\mu$ and $\nu$, respectively, $\tau$ is a transposition, and $b\tau=a$. Since $\tau^2=1\!\!1_{n^+,n^-}$, the equality $b\tau=a$ is equivalent to the equality $a\tau=b$. We conclude the relation
$$
w^+_{\mu\nu}\,|C_\mu|=w^+_{\nu\mu}\,|C_\nu|.
$$
This is equivalent to the self-adjointness of $W^+$:
$$(W^+C_\nu,C_\mu)=w^+_{\mu,\nu}n^+!\;n^-!\;|C_\mu|=w^+_{\nu,\mu}n^+!\;n^-!\;|C_\nu|=(W^+C_\mu,C_\nu),$$
and the proposition is proved.

\section*{Appendix}

The problem of enumerating simple real meromorphic functions on not necessarily separating real curves
is treated in~\cite{N10} under the name of ``computation of  the Hurwitz numbers of a disk''
(see also \cite{MMN13a,MMN13b,MMN14}).
The two cases, namely, those of simple real functions on arbitrary curves and on only separating
ones are parallel. In this appendix, we review the results of~\cite{N10}, but in
notation adopted to the content of the present paper, and we compare
the two lines paying special attention to the points where they differ.

To a real curve $(C,\tau)$ we associate the quotient $D=C/\tau$ which is a surface with boundary.
The surface $D$ is orientable if and only if the curve is separating. In the separating case the surface $D$ is also diffeomorphic to either of the two connected components into which the curve of real points
$C^\tau$ splits
 the surface~$C$. A real meromorphic function $f:(C,\tau)\to(\C P^1,\sigma)$ defines the corresponding map of quotient surfaces $C/\tau\to\C P^1/\sigma$, and this correspondence establishes a bijection between isomorphism classes of real functions and coverings of a disk in a sense of Dold--Smith treated in~\cite{N10}.

\begin{figure}
%\centerline{\ig[scale=.3]{IMG_6no.jpg}}
\setlength{\unitlength}{0.5mm}
\thicklines
\begin{picture}(400,110)(0,0)
\qbezier(0,80)(0,110)(30,110)
\qbezier(0,80)(0,60)(40,60)
\qbezier(40,60)(80,60)(80,80)
\qbezier(80,80)(80,110)(60,113)
\qbezier(60,113)(52,113)(47,110)
\qbezier(30,110)(60,110)(60,90)
\qbezier(60,90)(60,80)(40,80)
\qbezier(40,80)(20,80)(20,90)
\qbezier(20,90)(20,100)(39,107)
\qbezier(30,70)(30,75)(40,75)
\qbezier(40,75)(50,75)(50,70)
\qbezier(30,70)(30,65)(40,65)
\qbezier(40,65)(50,65)(50,70)
\put(40,65){\vector(0,-1){5}}
\put(40,80){\vector(0,-1){5}}
\put(67,64){\vector(-1,1){17}}
\put(57,83){\vector(1,-2){10}}
\put(40,60){\circle*{2}}
\put(40,65){\circle*{2}}
\put(40,75){\circle*{2}}
\put(40,80){\circle*{2}}
\put(57,83){\circle*{2}}
\put(67,64){\circle{2}}
\put(50,81){\circle*{2}}
\put(37,54){\mbox{\footnotesize $5$}}
\put(39,66){\mbox{\footnotesize $4$}}
\put(35,75){\mbox{\footnotesize $3$}}
\put(41,81){\mbox{\footnotesize $2$}}
\put(59,82){\mbox{\footnotesize $6$}}
\put(63,60){\mbox{\footnotesize $\infty$}}
\put(45,78){\mbox{\footnotesize $1$}}
\qbezier(100,80)(100,110)(130,110)
\qbezier(100,80)(100,60)(140,60)
\qbezier(140,60)(180,60)(180,80)
\qbezier(180,80)(180,110)(160,113)
\qbezier(160,113)(152,113)(147,110)
\qbezier(130,110)(160,110)(160,90)
\qbezier(160,90)(160,80)(140,80)
\qbezier(140,80)(120,80)(120,90)
\qbezier(120,90)(120,100)(139,107)
\qbezier(130,70)(130,75)(140,75)
\qbezier(140,75)(150,75)(150,70)
\qbezier(130,70)(130,65)(140,65)
\qbezier(140,65)(150,65)(150,70)
\qbezier(140,60)(160,70)(157,83)
\put(140,60){\vector(0,1){5}}
\put(140,75){\vector(0,1){5}}
\put(123,83){\vector(-1,-1){16}}
\put(144,62){\vector(-2,-1){4}}
\put(140,60){\circle{2}}
\put(140,65){\circle*{2}}
\put(140,75){\circle*{2}}
\put(140,80){\circle*{2}}
\put(123,83){\circle*{2}}
\put(107,66){\circle*{2}}
\put(157,83){\circle*{2}}
\put(137,54){\mbox{\footnotesize $\infty$}}
\put(139,66){\mbox{\footnotesize $1$}}
\put(135,75){\mbox{\footnotesize $2$}}
\put(141,81){\mbox{\footnotesize $3$}}
\put(115,79){\mbox{\footnotesize $4$}}
\put(103,59){\mbox{\footnotesize $5$}}
\put(159,79){\mbox{\footnotesize $6$}}
\qbezier(200,80)(200,110)(230,110)
\qbezier(200,80)(200,60)(240,60)
\qbezier(240,60)(280,60)(280,80)
\qbezier(280,80)(280,110)(260,113)
\qbezier(260,113)(252,113)(247,110)
\qbezier(230,110)(260,110)(260,90)
\qbezier(260,90)(260,80)(240,80)
\qbezier(240,80)(220,80)(220,90)
\qbezier(220,90)(220,100)(239,107)
\qbezier(230,70)(230,75)(240,75)
\qbezier(240,75)(250,75)(250,70)
\qbezier(230,70)(230,65)(240,65)
\qbezier(240,65)(250,65)(250,70)
\qbezier(240,60)(215,70)(230,81)
\put(240,65){\vector(0,-1){5}}
\put(240,80){\vector(0,-1){5}}
\put(223,83){\vector(-1,-1){16}}
\put(226,77){\vector(1,1){4}}
\put(240,60){\circle{2}}
\put(240,65){\circle*{2}}
\put(240,75){\circle*{2}}
\put(240,80){\circle*{2}}
\put(223,83){\circle*{2}}
\put(207,66){\circle*{2}}
\put(230,81){\circle*{2}}
\put(237,54){\mbox{\footnotesize $\infty$}}
\put(239,66){\mbox{\footnotesize $6$}}
\put(235,75){\mbox{\footnotesize $5$}}
\put(241,81){\mbox{\footnotesize $4$}}
\put(215,79){\mbox{\footnotesize $2$}}
\put(203,59){\mbox{\footnotesize $3$}}
\put(231,82){\mbox{\footnotesize $1$}}
\qbezier(140,30)(140,50)(180,50)
\qbezier(220,30)(220,50)(180,50)
\qbezier(220,30)(220,0)(200,0)
\qbezier(178,12)(190,0)(200,0)
\qbezier(140,30)(140,10)(160,10)
\qbezier(160,10)(177,10)(185,26)
\qbezier(186,23)(195,27)(207,27)
\qbezier(207,27)(213,27)(213,15)
\qbezier(213,15)(213,5)(200,5)
\qbezier(200,5)(190,5)(183,14)
\qbezier(174,18)(167,30)(160,30)
\qbezier(160,30)(153,30)(153,23)
\qbezier(153,23)(153,16)(160,16)
\qbezier(160,16)(170,16)(181,22)
\qbezier(185,26)(192,38)(180,38)
\qbezier(180,38)(175,38)(175,28)
\qbezier(175,28)(175,25)(177,21)
\put(160,30){\vector(1,1){20}}
\put(180,50){\vector(0,-1){12}}
\put(188,32){\vector(2,-1){13}}
\put(212,22){\vector(1,1){8}}
\put(180,50){\circle{2}}
\put(201,26){\circle*{2}}
\put(160,30){\circle*{2}}
\put(220,30){\circle*{2}}
\put(188,32){\circle*{2}}
\put(212,22){\circle*{2}}
\put(180,38){\circle*{2}}
\put(176,52){\mbox{\footnotesize $\infty$}}
\put(199,19){\mbox{\footnotesize $3$}}
\put(189,33){\mbox{\footnotesize $2$}}
\put(215,30){\mbox{\footnotesize $5$}}
\put(177,31){\mbox{\footnotesize $1$}}
\put(210,26){\mbox{\footnotesize $4$}}
\put(156,31){\mbox{\footnotesize $6$}}
\qbezier(40,30)(40,50)(80,50)
\qbezier(120,30)(120,50)(80,50)
\qbezier(120,30)(120,0)(100,0)
\qbezier(78,12)(90,0)(100,0)
\qbezier(40,30)(40,10)(60,10)
\qbezier(60,10)(77,10)(85,26)
\qbezier(86,23)(95,27)(107,27)
\qbezier(107,27)(113,27)(113,15)
\qbezier(113,15)(113,5)(100,5)
\qbezier(100,5)(90,5)(83,14)
\qbezier(74,18)(67,30)(60,30)
\qbezier(60,30)(53,30)(53,23)
\qbezier(53,23)(53,16)(60,16)
\qbezier(60,16)(70,16)(81,22)
\qbezier(85,26)(92,38)(80,38)
\qbezier(80,38)(75,38)(75,28)
\qbezier(75,28)(75,25)(77,21)
\put(80,50){\vector(-1,-1){20}}
\put(80,38){\vector(0,1){12}}
\put(101,26){\vector(-2,1){13}}
\put(120,30){\vector(-1,-1){8}}
\put(80,50){\circle{2}}
\put(101,26){\circle*{2}}
\put(60,30){\circle*{2}}
\put(120,30){\circle*{2}}
\put(88,32){\circle*{2}}
\put(112,22){\circle*{2}}
\put(80,38){\circle*{2}}
\put(76,52){\mbox{\footnotesize $\infty$}}
\put(99,19){\mbox{\footnotesize $4$}}
\put(89,33){\mbox{\footnotesize $5$}}
\put(115,30){\mbox{\footnotesize $2$}}
\put(77,31){\mbox{\footnotesize $6$}}
\put(110,26){\mbox{\footnotesize $3$}}
\put(56,31){\mbox{\footnotesize $1$}}
\end{picture}
\caption{All the~$5$ nonorientable diagrams of genus~$2$ with a single
boundary pole of degree~$3$; together with the orientable diagrams in Fig.~\ref{fd} they
form the~$9$ diagrams on both separating and nonseparating real curves}\label{fdnsep}
\end{figure}
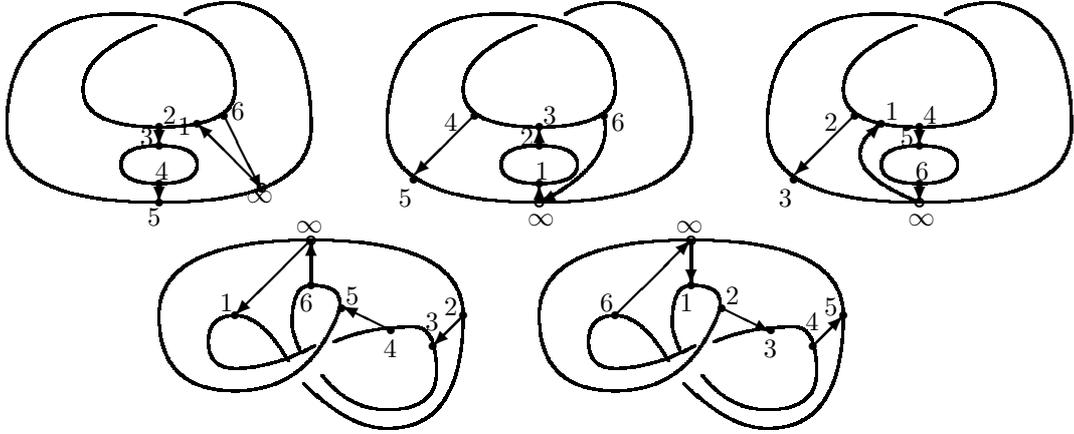

If a real function is simple, then it can be uniquely reconstructed from its diagram similar to the one discussed in Sect.~\ref{sDD}. The definition of the diagram of a simple real function on a not necessarily
separating real curve repeats that of Definition~\ref{defdiag} with the only difference:
the diagram itself is not assumed to be orientable. For example, for degree $3$ real functions on a genus~$3$ curve with a single pole we have in addition to oriented diagrams of Fig.~\ref{fd} five more diagrams,
 which are not orientable; they are depicted in Fig.~\ref{fdnsep}.
 The underlying surface of the upper three of the nonorientable diagrams is the M\"obius
  band, while for the two others it is the Klein bottle without a disk.

Similarly to the separating case, the poles of a simple real function split into real ones and pairs of conjugate non-real poles. However, since $C^\tau$ is nonorientable,
the \emph{sign of a real pole} is well defined for a pole of \emph{even order only}:
 it is positive (resp. negative) if the corresponding critical point is a local minimum (resp., local maximum) of the restriction of a function to the set of real points.
 Thus, the \emph{ramification type} of $f$ at infinity is a quadruple $\mu=(\kappa^+,\kappa^-,\kappa,\lambda)$ where $\kappa^+$ and $\kappa^-$ are the partitions formed by the even parts corresponding to the orders of positive and negative real poles,
 respectively, $\kappa$ is the partition formed by odd parts corresponding to the orders of poles of odd orders, and $\lambda$ is a partition formed by the orders of pairs of conjugate non-real poles.
 The corresponding simple purely real Hurwitz numbers that we denote by $\tilde h^{\R}_{m;\mu}$, $\tilde h^{\R\circ}_{m;\mu}$, and the corresponding generating functions denoted by $\tilde H^{\R}$ and $\tilde H^{\R\circ}=\exp(\tilde H^{\R})$ are defined in a way similar to the separating case of Sect.~\ref{sFSHN}.
 But since the signs of poles of odd orders are not distinguished, the set of independent variables entering the functions $\tilde H^{\R}$ and  $\tilde H^{\R\circ}$ consists of the following ones: $p_{2i}^+$, $p_{2i}^-$, $p_{2i-1}$, $q_i$, $i=1,2,3,\dots$ (in~\cite{N10} these variables are denoted by $\acute{p}_i$, $\grave{p}_i$, $\bar{p}_i$, $\dot{p}_i$, respectively).

\begin{theorem}[\cite{N10}]\label{tcjdpm}
The generating function $\tilde H^{\R\circ}$ satisfies the differential equation
\begin{equation}\label{ecjd}
\frac{\p \tilde H^{\R\circ}}{\partial u}=\tilde W(\tilde H^{\R\circ});
\end{equation}
here the linear partial differential operator~$\tilde W$ is defined as follows:
\begin{eqnarray*}
\tilde W&=&\sum _{i,j} \left(2 p_{2 i+2 j-1} \hbar  \frac{\partial ^2}{\partial p_{2 i-1}\, \partial p_{2 j}^+}+
\frac{1}{2} p_{2 i+2 j-2}^- \hbar  \frac{\partial ^2}{\partial p_{2 i-1}\, \partial p_{2 j-1}}+
2 p_{2 i+2 j}^+ \hbar  \frac{\partial ^2}{\partial p_{2 i}^+\, \partial p_{2 j}^+}\right)\\
&&+
\sum _{i,j}\left( p_{2 i-1} p_{2 j-1} \frac{\partial }{\partial p_{2 i+2 j-2}^-}+p_{2 i-1} p_{2 j}^+ \frac{\partial }{\partial p_{2 i+2 j-1}}+p_{2 i}^+ p_{2 j}^+ \frac{\partial }{\partial p_{2 i+2 j}^+}\right)\\
&&+
\sum_{i=1}^\infty\left(i\hbar^{-2} p^+_{2i}\frac{\p}{\p q_i}+q_i\frac{\p}{\p p^+_{2i}}\right).
\end{eqnarray*}
\end{theorem}

In contrast to the operators $W^+,W^-,W$, the operator $\tilde W$ preserves only grading rather than bigrading.

The proof in~\cite{N10} uses an algebraic model for real Hurwitz numbers introduced in~\cite{AN06}, Sec.~5.
 This model is similar to that of Sec.~\ref{sTFT} (and even, in a sense, a bit simpler).
 Let us recall it here. Fix a finite set $N$ consisting of $|N|=n$ elements. A \emph{state} on $N$ is defined as an arbitrary involutive permutation, that is, a splitting of $N$ into a disjoint union of one and two-element subsets.
 A \emph{transition} is an ordered pair of states. For example, a \emph{transposition} is a
 transition whose states differ exactly by two elements that form a two-element subset
 in one of the states and are individual elements in the other state.
 Comparing with the definitions of Sect.~\ref{sTFT} we see that here
 the elements in~$N$ do not have signs. A transition can be depicted
 by means of a diagram similar to that in Fig.~\ref{ftrans}, with all signs dropped.

A \emph{transition type} is an orbit of the permutation group $S(n)$ acting on the set of transitions
by relabeling the elements of~$N$. One observes that the transition types
are in one-to-one correspondence with the ramification types of degree~$n$ real functions.
Transitions span an algebra denoted by $\tilde T_n$.
Denote also by $\tilde A_n=\tilde T_n^{S(n)}$ its $S(n)$-invariant subalgebra.
It is generated by the elements $C_\mu$ given as the sum of all transitions of a given type $\mu$.
An argument similar to that in Sect.~\ref{sTFT} shows that the Hurwitz number
$\tilde h^{\R\circ}_{m,\mu}$ is equal to the suitably rescaled coefficient of
$C_\mu$ in the product $C_2^m$ of transpositions. Thus, the operator $\tilde W$ of
Theorem~\ref{ecjd} describes the action of multiplication by $C_2$
in the algebra $\tilde A_n$, where $\tilde A_n$ is identified
with the space of weighted degree~$n$ monomials in the variables $p^+_{2i},p^-_{2i},p_{2i-1},q_{i}$.

\end{document}